\documentclass[11pt]{amsart}
\usepackage{amsxtra,tikz-cd,color}
\usepackage{amssymb}
\usepackage{mathtools}
\theoremstyle{plain}
\usepackage{underscore}

\newcommand{\cleqn}{\setcounter{equation}{0}}
\newcommand{\clth}{\setcounter{theorem}{0}}
\newcommand {\sectionnew}[1]{\section{#1}\cleqn\clth}
\newcommand{\nn}{\hfill\nonumber}
\newtheorem{theorem}{Theorem}[section]
\newtheorem{lemma}[theorem]{Lemma}
\newtheorem{definition-theorem}[theorem]{Definition-Theorem}
\newtheorem{proposition}[theorem]{Proposition}
\newtheorem{corollary}[theorem]{Corollary}
\newtheorem{definition}[theorem]{Definition}
\newtheorem{example}[theorem]{Example}
\newtheorem{remark}[theorem]{Remark}
\newtheorem{conjecture}[theorem]{Conjecture}

\newcommand \bth[1] { \begin{theorem}\label{t#1} }
\newcommand \ble[1] { \begin{lemma}\label{l#1} }

\newcommand \bpr[1] { \begin{proposition}\label{p#1} }
\newcommand \bco[1] { \begin{corollary}\label{c#1} }
\newcommand \bde[1] { \begin{definition}\label{d#1}\rm }
\newcommand \bex[1] { \begin{example}\label{e#1}\rm }
\newcommand \bre[1] { \begin{remark}\label{r#1}\rm }
\newcommand \bcj[1] { \begin{conjecture}\label{j#1}\rm }

\renewcommand {\eth} { \end{theorem} }
\newcommand {\ele} { \end{lemma} }

\newcommand {\epr} { \end{proposition} }
\newcommand {\eco} { \end{corollary} }
\newcommand {\ede} { \end{definition} }
\newcommand {\eex} { \end{example} }
\newcommand {\ere} { \end{remark} }
\newcommand {\ecj} { \end{conjecture} }

\newcommand {\enota} { \end{notation} }
\newcommand \thref[1]{Theorem \ref{t#1}}
\newcommand \leref[1]{Lemma \ref{l#1}}
\newcommand \prref[1]{Proposition \ref{p#1}}

\newcommand \deref[1]{Definition \ref{d#1}}
\newcommand \exref[1]{Example \ref{e#1}}

\def \Cset {{\mathbb C}}

\def \Zset {{\mathbb Z}}

\def \Qset {{\mathbb Q}}

\def \CC {\mathcal{C}}

\def \VV {{\mathcal{V}}}
\def \UU {{\mathcal{U}}}

\def \OO {{\mathcal{O}}}

\def \al {\alpha}

\def \la {\lambda}

\def \Om {\Omega}
\def \ga {\gamma}

\def \vp {\varphi}
\def \ep {\varepsilon}

\def \mt  {\mapsto}

\def \ol {\overline}
\def \wt {\widetilde}

\def \id { {\mathrm{id}} }
\def \Stab { {\mathrm{Stab}} }

\def \Lie { {\mathrm{Lie \,}} }

\def \g  {\mathfrak{g}}   

\def \mm  {\mathfrak{m}}

\newcommand{\kk}{\Bbbk}

\newcommand{\ee}{\epsilon}

\DeclareMathOperator \Aut { {\mathrm{Aut}} }

\DeclareMathOperator \rad { {\mathrm{rad}}}

\DeclareMathOperator \charr { {\mathrm{char}} }

\DeclareMathOperator \opp { {\mathrm{op}} }

\DeclareMathOperator \MaxSpec { {\mathrm{MaxSpec}}}
\DeclareMathOperator \ord { {\mathrm{ord}} }

\DeclareMathOperator \diag { {\mathrm{diag}} }
\DeclareMathOperator \Ker { {\mathrm{Ker}} }
\DeclareMathOperator \coker { {\mathrm{coker}} }
\DeclareMathOperator \tr { {\mathrm{tr}} }

\DeclareMathOperator \End { {\mathrm{End}} }

\DeclareMathOperator \Hom { {\mathrm{Hom}} }

\DeclareMathOperator \Ann { {\mathrm{Ann}} }

\DeclareMathOperator \Irr { {\mathrm{Irr}} }

\DeclareMathOperator \Ext { {\mathrm{Ext}} }
\DeclareMathOperator \alg {{\mathrm{alg}}}

\DeclareMathOperator \reg  { {\mathrm{reg}}}
\DeclareMathOperator \Sd {{\mathrm{Sd}}}

\DeclareMathOperator \modd {{\mathrm{mod}}}

\begin{document}
\title
{The lowest discriminant ideal of a Cayley--Hamilton Hopf algebra}
\author[Zhongkai Mi]{Zhongkai Mi}
\address{Shanghai Center for Mathematical Sciences \\
Fudan University \\
Shanghai, 200438 \\
People's Republic of China}
\email{zhongkai\_mi@fudan.edu.cn}
\author[Quanshui Wu]{Quanshui Wu}
\thanks{The research of Q.W. has been supported by the NSFC (Grant No. 12471032) and the National Key Research and Development Program of China (Grant No. 2020YFA0713200). The research of Z.M. and M.Y. has been supported by NSF grants DMS-2131243 and DMS--2200762.}
\address{School of Mathematical Sciences \\
Fudan University \\
Shanghai, 200433 \\
People's Republic of China}
\email{qswu@fudan.edu.cn}
\author[Milen Yakimov]{Milen Yakimov}
\address{
Department of Mathematics \\
Northeastern University \\
Boston, MA 02115 \\
U.S.A.}
\email{m.yakimov@northeastern.edu}
\date{}
\begin{abstract} Discriminant ideals of noncommutative algebras $A$, which are module finite over a central subalgebra $C$,
are key invariants that carry important information about $A$, such as the sum of the squares of the dimensions of its irreducible modules with a given
central character. There has been substantial research on the computation of discriminants, but
very little is known about the computation of discriminant ideals. In this paper we carry out a detailed investigation of the
lowest discriminant ideals of Cayley--Hamilton Hopf algebras in the sense of De Concini, Reshetikhin, Rosso and Procesi,
whose identity fiber algebras are basic. The lowest discriminant ideals are the most complicated ones, because they
capture the most degenerate behaviour of the fibers in the exact opposite spectrum of the picture from the Azumaya locus.
We provide a description of the zero sets of the lowest discriminant ideals of Cayley--Hamilton Hopf algebras in terms of maximally stable
modules of Hopf algebras, irreducible modules that are stable under tensoring with the maximal possible number
of irreducible modules with trivial central character. In important situations, this is shown to be governed by the actions of the
winding automorphism groups. The results are illustrated with applications to the group algebras of central
extensions of abelian groups, big quantum Borel subalgebras at roots of unity and quantum coordinate rings at roots of unity.
\end{abstract}
\maketitle
\sectionnew{Introduction}
\subsection{Discriminants and discriminant ideals of noncommutative algebras}
In many important situations in representation theory, commutative and noncommutative algebra, one is lead to consider
algebras $A$ over a commutative ring $\kk$ that are module finite over a Noetherian central subalgebra $C$. Following
\cite[Definition B.1.5]{B}, we call them {\em{Noether algebras}}.
By various general constructions, such algebras are equipped with
a {\em{trace function}} $\tr : A \to C$, which is $C$-linear and cyclic:
\[
\tr(c a) = c \tr(a), \quad \tr(ab) = \tr (ba), \quad \forall \, a, b \in A, \; c \in C.
\]
In this setting, one considers the {\em{discriminant ideals}} \cite{Rei}
\[
D_n(A/C, \tr) := \big( \det [ \tr(a_i a_j) ]_{i,j =1}^n \mid (a_1, \ldots, a_n) \in A^n \big)
\]
and the {\em{modified discriminant ideals}} \cite{CPWZ}
\[
MD_n(A/C, \tr) := \big( \det [ \tr(a_i a'_j) ]_{i,j =1}^n \mid (a_1, \ldots, a_n), (a'_1, \ldots, a'_n) \in A^n \big),
\]
where as usual, for a subset $S \subseteq C$,  $(S)$ denotes the ideal of $C$ generated by $S$. 
Both types of ideals are ideals of the central subalgebra $C$ of $A$. 
If $A$ is a free $C$-module of rank $N$, then one also considers the {\em{discriminant}} of $A$ over $C$ given by \cite{Rei}
\[
\Delta(A/C, \tr) = \det [\tr(b_i b_j) ]_{i,j =1}^N,
\]
where $\{b_1, \ldots, b_N \}$ is a $C$-basis of $A$. The discriminant is uniquely defined
up to a multiplication by the square of a unit of $C$ and, in this special case,
\begin{equation}
\label{relation}
D_N(A/C, \tr) = MD_N(A/C, \tr) = (\Delta(A/C, \tr)).
\end{equation}
Discriminants were first defined by Dedekind in the case when $A$ is the of ring integers of an algebraic
number field and $C = \Zset$, and since, have played a key role in algebraic number theory.
In noncommutative algebra, discriminants play a major role in the theory of orders \cite{Rei}, the classification of automorphism
groups of algebras \cite{CPWZ}, and the Zariski cancellation problem \cite{BZ1}. A deep relation between the
zero sets of the discriminant (modified discriminant ideals) of an algebra and the dimensions of its
irreducible modules was found in \cite{BY}.

Discriminants of noncommutative algebras have been intensely studied in recent years and many methods have been developed for their
evaluation: smash product techniques \cite{GKM}, Poisson geometry \cite{NTY1,LeYa}, cluster algebras \cite{NTY2}, reflexive hulls
\cite{CGWZ} and others.

However, extremely little is presently known about discriminant ideals. They play a much more important role than discriminants and
are defined in much wider generality without a freeness assumption on the noncommutative algebra $A$ over $C$. 
This assumption is in fact very often violated. For example, root of unity quantum cluster algebras have natural trace functions 
to canonical central subalgebras \cite[Theorem B]{HLY}, and freeness only holds for special families of cluster algebras.
Discriminant ideals even carry more information than discriminants in the case when
the latter are defined, because the relation \eqref{relation} only concerns the top discriminant ideal, while the other
ideals
\[
\{ D_k(A/C, \tr), MD_k(A/C, \tr) \mid 1 \leq k < N \}
\]
carry information not seen by the discriminant $\Delta(A/C, \tr)$.

In this paper we investigate the zero set of the {\em{lowest discriminant ideal}} of a Noether algebra with trace $(A,C, \tr)$.
For all $\mm \in \MaxSpec C$ (the maximum spectrum of $C$), consider the {\em{fiber algebras}}
\[
A/ \mm A.
\]
They are finite dimensional by the assumption that $A$ is module finite over $C$. Denote by
\[
\Irr(A/ \mm A)
\]
the {\em{isomorphism classes of irreducible modules}} of $A / \mm A$. If we are in one of the following two general situations
\begin{enumerate}
\item[(1)] $A$ is a prime affine algebra of PI degree $n$ over an algebraically closed field, which is module over its center $C:=Z(A)$ 
and coincides with its trace ring (e.g. $A$ is a {\em{maximal order}})  or
\item[(2)] $(A,C, \tr)$ is a finitely generated {\em{Cayley--Hamilton algebra of degree $n$}} in the sense of Procesi \cite{P,DP}, see
Section \ref{CH} for details,
\end{enumerate}
and $\charr \kk \notin [1,n]$, then by \cite[Main Theorem (e) and Theorem 4.1(b)]{BY},
\begin{align}
\label{VVk}
\VV_k \; := \; \VV(D_k(A/C, \tr)) &= \VV(MD_k(A/C, \tr))
\\
&= \Big\{ \mm \in \MaxSpec C \mid \sum_{V \in \Irr(A/ \mm A)} \dim(V)^2 < k \Big\}
\nn
\end{align}
for all positive integers $k$. All important families of Noether algebras $(A,C)$ fall in one of the two classes (1)--(2).

We have
\[
MD_k(A/C, \tr) \supseteq MD_{k+1}(A/C, \tr) \quad \mbox{for all} \quad k \geq 1, 
\]
and consequently, 
\[
\VV_1 \subseteq \VV_2 \subseteq \ldots \subseteq \VV_N \subsetneq \VV_{N+1} = \MaxSpec C
\]
for some positive $N$. It was proved in \cite{BY} that in situation (1), $N=n^2$, where $n$ is the PI degree of $A$, and that the complement
\[
\MaxSpec C \backslash \VV_N
\]
of the zero set of the highest discriminant ideal $D_N(A/C, \tr)$ is precisely the {\em{Azumaya locus}}
of $A$, \cite[Main Theorem (a)]{BY}. A similar statement holds in situation (2),
\cite[Theorem 4.1(a)]{BY}.

Denote the {\em{lowest modified discriminant ideal}} $MD_\ell(A/C, \tr)$ to be the one, corresponding to the unique positive
integer $\ell$ for which
\[
\varnothing = \VV_1 = \VV_2 = \ldots =  \VV_{\ell-1} \subsetneq \VV_\ell.
\]
By \eqref{VVk}, the same integer captures the
{\em{lowest discriminant ideal}} $D_\ell(A/C, \tr)$.

The goal of the paper is to address the following fundamental problem:
\medskip

\noindent
{\bf{Main Problem.}} {\em{Describe the zero set of the lowest discriminant ideal}}
\[
\VV_\ell \subset \MaxSpec C
\]
{\em{of a Cayley--Hamilton algebra $(A,C, \tr)$.}}
\medskip

We start with noting that the lowest discriminant ideal is much more complicated than the highest discriminant ideal. While
the complement of the zero set of the latter captures the generic fibers of the algebra $A$ for which $A / \mm A$ is a simple algebra
(or more generally a semisimple algebra, depending on the choice of $C$),
the zero set of the minimal discriminant ideal consists of the most degenerate fibers. These are precisely the fibers of
highest interest in representation theory; in special cases the most degenerate fiber algebras $A/ \mm A$
are isomorphic to restricted universal enveloping algebras in positive characteristic and small quantum groups at roots of unity.

By a {\em{character}} of a $\kk$-algebra $A$ we will mean a $\kk$-algebra homomorphism $A \to \kk$.
The characters of $A$
are in bijection with the isomorphism classes of 1-dimensional modules of $A$. We will identify the two sets and, for a
character $\phi : A \to \kk$, we will denote by the same letter the corresponding 1-dimensional module of $A$.
Given a character $\phi : A \to \kk$, denote the maximal ideal
\[
\mm_\phi := \Ker \phi \; \; \mbox{of} \; \; A.
\]

For a Noether Hopf algebra $H$ with respect to a central Hopf subalgebra $C \subseteq H$, we denote by $\ep: H \to \kk$ and $\ol{\ep} := \ep|_C : C \to \kk$
the counits of $H$ and $C$, respectively.
{\em{We view such a Hopf algebra $H$ as the total space of the finite dimensional deformations $H/ \mm H$ {\em{(}}for $\mm \in \MaxSpec C${\em{)}}
of the finite dimensional Hopf algebra $H/\mm_{\ol{\ep}}$, see Section \ref{act}.}}
\subsection{Statements of the main results in the paper}
We address the Main Problem under the following general assumptions:
\begin{enumerate}
\item[(1)] $(H, C, \tr)$ is a finitely generated Cayley--Hamilton Hopf algebra in the sense of De Concini, Procesi, Reshetikhin and Rosso
\cite{DPRR}, which means that $(H, C, \tr)$ is a Cayley--Hamilton algebra and $C$ is a central Hopf subalgebra of $H$, see \deref{CHHopfalg}, and
\item[(2)] the identity fiber $H/ \mm_{\ol{\ep}}H$ is a basic algebra.
\end{enumerate}

By \cite[Theorem 4.5]{DP}, assumption (1) implies that $H$ is module finite over $C$ and $C$ is a finitely generated $\kk$-algebra. So, in particular,
$(H,C)$ is a Noether algebra.

In practice, assumption (1) encompasses the majority of PI Hopf algebras that play a role in representation theory, mathematical physics and topology.
One of the main directions in the area of Hopf algebras is the classification of pointed finite dimensional Hopf algebras \cite{AS}. 
It resulted in the construction of many remarkable classes
of such algebras, see \cite{AA,HS}. The assumption (2) is equivalent to saying
that we work in the situation that the dual of the identity fiber algebra $(H/ \mm_{\ol{\ep}})^*$ is a pointed Hopf algebra. In other words,
{\em{we consider the discriminant ideals of all Hopf algebras $H$ which are the total spaces of the deformations of
the dual of any finite dimensional pointed Hopf algebra $(H/ \mm_{\ol{\ep}})^*$}}.

Our first theorem addresses the action of the isomorphism classes of the irreducible modules of the identity fiber algebra $H/ \mm_{\ol{\ep}} H$
on the isomorphism classes of the irreducible modules of the fiber algebras $H/ \mm H$ and, in particular, it describes the critical situation for its stabilizers.
\medskip

\noindent
{\bf{Theorem A.}} {\em{Assume that $H$ is a finitely generated Hopf algebra over an algebraically closed field $\kk$
and $C$ is a central Hopf subalgebra such that $H$ is a finitely generated $C$-module. Then the following hold:}}
\begin{enumerate}
\item[(a)] {\em{For all $V, W \in \Irr(H/\mm H)$, there exist $M', M'' \in \Irr(H/\mm_{\ol{\ep}} H)$ such that $W$ is a quotient of $M' \otimes V$
and a submodule of $M'' \otimes V$.}}
\item[(b)] {\em{The group
\begin{equation}
\label{G0}
G_0:=G( (H/\mm_{\ol{\ep}}H)^\circ ) \subseteq \Irr (H/\mm_{\ol{\ep}}H)
\end{equation}
acts on $\Irr(H/\mm H)$ by $V \in \Irr(H/\mm H) \mt \chi \otimes V$ for $\chi \in G_0$.}}
\item[(c)] {\em{Every module $V \in \Irr(H/ \mm H)$ for $\mm \in \MaxSpec C$ satisfies}}
\begin{equation}
\label{Stab-inequal}
|\Stab_{G_0}(V)| \leq \dim(V)^2.
\end{equation}
\item[(d)]
{\em{An equality in \eqref{Stab-inequal} holds {\em{(}}we call such a module {\bf{maximally stable}}{\em{)}} if and only if
one of the following two equivalent conditions holds:}}
\begin{enumerate}
\item[(ii)] {\em{$V \otimes V^*$ is a direct sum of nonisomorphic 1-dimensional $H$-modules.}}
\item[(iii)] {\em{$V \otimes V^* \cong \bigoplus_{\chi \in \Stab_{G_0}(V)} \chi$.}}
\end{enumerate}
\item[(e)] {\em{If $V \in  \Irr(H/\mm H)$ is a maximally stable irreducible module, then the primitive quotient
\[
H/ \Ann_H(V)
\]
is isomorphic to a twisted group algebra
\[
\kk_{\gamma_V} \Stab_{G_0}(V)
\]
for a canonically defined 2-cocycle $\gamma_V : \Stab_{G_0}(V) \times \Stab_{G_0}(V) \to \kk^*$,
see \thref{twisted-rep}{\em{(}}a{\em{)}}.
Both algebras are isomorphic to $\End_\kk(V)$.}}
\end{enumerate}
{\em{If, in addition, the identity fiber algebra $H/ \mm_{\ol{\ep}} H$ is a basic algebra, then we have:}}
\begin{enumerate}
\item[(f)] \cite[Proposition III.4.11]{BG-book}
{\em{The inclusion in part}} ({\em{b}}) {\em{is an equality and the action in part (b) is transitive.}}
\end{enumerate}
\medskip

\noindent
{\bf{Remark.}} Recall that a {\em{group of central type}} \cite{DJ,HI,LY} is a finite group $G$ which has an irreducible complex character $\varphi$ 
with degree $\varphi(1) = \sqrt{[G:Z(G)]}$ (as usual, $Z(G)$ denotes the center of $G$ and $[G:Z(G)]$ 
its index). For such a group, $G/Z(G)$ has a 2-cocycle $\gamma$ such that the twisted group algebra $\kk_\gamma (G/Z(G))$
is simple, and vice versa, if $H$ is a finite group possessing a 2-cocycle for which the twisted group ring 
$\kk_\gamma H$ is simple, then $H \cong G/Z(G)$ for a group $G$ of central type \cite[Theorem 1]{DJ}.

A conjecture of Iwahori and Matsumoto that a group of central type is solvable. The conjecture was proved (by using the classification of finite groups) 
by Liebler--Yellen \cite{LY} and Howlett--Isaacs \cite{HI} (the first proof contained a gap which was repaired in the second). By these results, 
in charateristic 0, the groups $\Stab_{G_0}(V)$ that appear in Theorem A(e) are all solvable. 
\medskip

The next is our first theorem on the zero sets of the lowest discriminant ideals of Cayley--Hamilton Hopf algebras
with basic identity fiber.
\medskip

\noindent
{\bf{Theorem B.}}
Assume that $(H,C,\tr)$ is a Cayley--Hamilton Hopf algebra of degree $n$
over an algebraically closed field $\kk$ of characteristic $\charr \kk \notin [1, n]$,
such that $H$ is a finitely generated $\kk$-algebra and
the identity fiber $H/ \mm_{\ol{\ep}} H$ is a basic algebra. The following hold:
\begin{enumerate}
\item[(a)] For all $\mm \in \MaxSpec C$,
\[
\min \big\{ k \in \Zset_+ \mid D_k(H/C, \tr)(\mm) = 0 \big\} = \frac{|G_0| \dim(V)^2}{|\Stab_{G_0}(V)|} +1
\]
for any $V \in \Irr(H/ \mm H)$, recall the definition \eqref{G0} of the group $G_0$.
\item[(b)] The lowest discriminant ideal of $(H,C, \tr)$ is of level
\[
\ell = |G_0| +1.
\]
\item[(c)] The following are equivalent for $\mm \in \MaxSpec C$:
\begin{enumerate}
\item[(i)] $\mm$ belongs to the zero set $\VV_{|G_0| +1}$ of the lowest discriminant ideal of $(H,C, \tr)$;
\item[(ii)] There exists $V \in \Irr(H/ \mm H)$ that is maximally stable;
\item[(iii)] All modules $V \in \Irr(H/ \mm H)$ are maximally stable.
\end{enumerate}
\item[(d)] Under the assumptions of the theorem, but 
without the assumption that the identity fiber $H/ \mm_{\ol{\ep}} H$ is a basic algebra, we have that
the group $G_0$ contains a copy of the cyclic group of order equal to the integral order $io(H)$ of $H$ and
for the order $\ell$ of the lowest discriminant ideal of $H$, we have
\[
\ell > io(H).
\]
\end{enumerate}
\medskip

By the above remark, in characteristic 0, to every point of the zero set $\VV_{|G_0| +1}$ of the lowest discriminant ideal of $(H,C, \tr)$, one can associate a canonical 
group of central type.  

By \cite[Theorems 0.1 and 0.2(1)]{WZ}, in the setting of Theorem B, $H$ is an {\em{AS-Gorenstein algebra}}, see \deref{AS-Gorenstein}. In this setting,
one defines \cite{LWZ} the left and right homological integrals of $H$, $\int^l$ and $\int^r$, which are certain $H$-bimodules,
and the integral order $io(H)$ which equals the minimal positive integer $n$ such that $(\int^r)^{\otimes n}$ is isomorphic to the trivial
trivial $H$-bimodule. This is the integer that appears in Theorem C(d), see \deref{hom-int} for details.
If $H$ is prime and regular, and has Gelfand--Kirillov dimension 1, then $io(H)$ equals the PI degree of $H$ by \cite[Theorem 0.2(a)]{LWZ}.

The estimate in Theorem B(d) is optimal in the sense that $\ell=|G_0|+1$ can be equal to $io(H)+1$. This is the case for \cite[Example 2.7]{LWZ} if we take $m=1$, $t=0$ and $C$ to be the center $\kk[x^n]$. In general, $io(H)$ can be much smaller than $|G_0|$. In fact, for the Hopf algebras $H$ in the examples in Sections \ref{subsec:fin-group}--\ref{subsec:inf-group}, 
$io(H)=1$.

Our second theorem on the zero set of the lowest discriminant ideal of a Cayley--Hamilton Hopf algebra $(H,C,\tr)$
proves that it contains an explicitly constructed subset given in terms of an automorphism group orbit and establishes that the
two sets coincide if every maximally stable irreducible $H$-module is 1-dimensional.
\medskip

\noindent
{\bf{Theorem C.}}
{\em{Assume that $(H,C,\tr)$ is a Cayley--Hamilton Hopf algebra satisfying the assumptions of Theorem B.}}
\begin{enumerate}
\item[(a)] {\em{For $\mm \in \MaxSpec C$ the following are equivalent:}}
\begin{enumerate}
\item[(i)] {\em{The algebra $H/\mm H$ has a 1-dimensional module;}}
\item[(ii)] {\em{The algebra $H/\mm H$ is basic;}}
\item[(iii)] {\em{$H/\mm H$ and $H/\mm_{\ol{\ep}} H$ are isomorphic as $\kk$-algebras;}}
\item[(iv)] {\em{$\mm$ belongs to the orbit of $\mm_{\ol{\ep}}$ under the group $\Aut_{\kk-\alg}(H,C)$
of all $\kk$-algebra automorphisms of $H$ that preserve $C$;}}
\item[(v)] {\em{$\mm$ belongs to the orbit of $\mm_{\ol{\ep}}$ under the left winding automorphism group $W_l(G(H^\circ))$;}}
\item[(vi)] {\em{$\mm$ belongs to the orbit of $\mm_{\ol{\ep}}$ under the right winding automorphism group $W_r(G(H^\circ))$.}}
\end{enumerate}
\item[(b)] {\em{If $\mm \in \MaxSpec C$ satisfies any of the six equivalent conditions in part {\em{(}}a{\em{)}}, then
it belongs to the zero set $\VV_{|G_0| +1}$ of the lowest discriminant ideal of $(H,C, \tr)$.}}
\item[(c)] {\em{If every maximally stable irreducible $H$-module is 1-dimensional, then the maximal ideal $\mm \in \MaxSpec C$
belongs to the zero set $\VV_{|G_0| +1}$ of the lowest discriminant ideal of $(H,C, \tr)$
if and only if it satisfies any of the six equivalents conditions in part {\em{(}}a{\em{)}}.}}
\end{enumerate}
\medskip

We refer the reader to Section \ref{sec:aut} for background on winding automorphisms of Hopf algebras.

We prove stronger results than those in Theorems B and C, which deal with the lowest level set of the square dimension
function
\[
\mm \mt \sum_{V \in \Irr(A/ \mm A)} \dim(V)^2
\]
for an arbitrary pair $(H,C)$ of a Hopf algebra $H$ and a central Hopf subalgebra $C$ such that
$H$ is module finite over $C$. Those results are obtained in Theorems \ref{tfirst-lowest-sd-set} and \ref{tsecond-lowest-Sd-set}.

In Section \ref{sec:examples} we present applications of Theorems A--C to group algebras of central extensions of abelian groups
and (big) quantum groups at roots of unity.
These applications answer the following naturally arising questions from these theorems:
\begin{enumerate}
\item[(Q1)] {\em{Is there a restriction on the dimensions of the maximally stable representations of a Cayley--Hamilton Hopf algebra?}}

For each maximally stable module $V$ of $H$, $\dim(V)^2$ equals $\Stab_{G_0}(V)$, and thus, $\dim(V)^2$ divides $|G_0|$.
We show that in general a Cayley--Hamilton Hopf algebra $H$ can have maximally stable representations $V$
of any dimension $\dim(V)$ such that
\[
\dim(V)^2 \mid |G_0|.
\]
\item[(Q2)] {\em{Is there a restriction on the number of maximally stable irreducible modules of a Cayley--Hamilton Hopf algebra $H$,
compared to all irreducible modules of $H$?}}

We show that the answer is no, and even more, there are Cayley--Hamilton Hopf algebras with the property that
all of their irreducible modules are maximally stable and at the same time the majority of those modules have dimensions
$>1$.

\item[(Q3)] {\em{Is the inclusion in Theorem C{\em{(}}b{\em{)}} proper in general?}}

We show that the answer is yes and that it can even happen that
\[
\Aut(H,C) . \mm_{\ol{\ep}} = \{ \mm_{\ol{\ep}} \},
\]
while the zero set of the lowest discriminant ideal of $(H,C, \tr)$ is the whole $\MaxSpec C$.

However, we also prove that in the important cases of the De Concini--Kac--Procesi big quantum Borel subgroups at roots of unity
\cite{DKP-solv,DP} and quantum function algebras at roots of unity \cite{DL,BG0}, there are no maximally stable representations
of dimension higher than 1, and thus, the inclusion in Theorem C(b) is an equality by Theorem C(c).
\end{enumerate}

The paper is organized as follows. Section \ref{background} contains background material.
In Section \ref{sec:act}, we study an action of the isomorphism classes of the irreducible modules 
of the identity fiber algebra $H/ \mm_{\ol{\ep}} H$ on the isomorphism classes of the irreducible modules 
of the fiber algebras $H/ \mm H$, define and investigate maximally stable modules, and prove Theorem A.
In Section \ref{sec:Sd}, we prove two characterizations of the lowest level set of the square dimension function of a Hopf
algebra $H$ which is module finite over a central Hopf subalgebra $C$, and used them to prove Theorems B and C.
In Section \ref{sec:examples} we present the above mentioned applications
\medskip

\noindent
{\bf{Notation.}} For a $\kk$-algebra $A$, we will denote by $\modd(A)$ the category of its finite dimensional (left) modules
and by $\Irr(A)$ the isomorphism classes of its irreducible modules. The center of $A$ will be denoted by $Z(A)$.
The algebra of $n \times n$ matrices with entries in a commutative ring $C$ will be denoted by $M_n(C)$.
The cardinality of a finite set $X$ will be denoted by $|X|$. For a finite abelian group $\Lambda$ and an algebraically closed field $\kk$,
the dual group of $\Lambda$ will be denoted by
\[
\widehat{\Lambda} : = \Hom (\Lambda, \kk^*).
\]

\noindent
{\bf{Acknowledgements.}} We are grateful to the anonymous referee for many remarks and suggestions which helped 
us to improve the exposition. 
\sectionnew{Background}
\label{background}
This section contains background material on the irreducible modules of algebras that are module finite over central subalgebras,
Cayley--Hamilton algebras, and winding automorphisms of Hopf algebras.
\subsection{Noether algebras}
\label{Noether-alg}
Assume that {\em{$A$ is a finitely generated algebra over an algebraically closed field $\kk$ and
$C$ is a central subalgebra such that $A$ is a finitely generated module over $C$, generated by $d$ elements}}.
The term {\em{Noether algebra}} was introduced for these algebras $A$ in \cite[Definition B.1.5]{B}.
(If $A$ is prime, then $A$ is an order in a central simple algebra by a special case of the Posner's Theorem
\cite[Section I.13.3]{BG-book}.)

By the Artin--Tate Lemma (see e.g. \cite[Section I.13.4]{BG-book}),
\begin{enumerate}
\item[(1)] $C$ is also finitely generated $\kk$-algebra, and hence is Noetherian.
\end{enumerate}
Kaplansky's Theorem  \cite[Section I.13.3]{BG-book} implies that
\begin{enumerate}
\item[(2)] every irreducible $A$-module $V$ is finite dimensional over $\kk$ and $\Ann_C(V) \in \MaxSpec C$, i.e., $V$ descends to
an irreducible $A/ \mm A$-module for some $\mm \in \MaxSpec C$,
\end{enumerate}
see \cite[Proposition III.1.1(4)]{BG-book}. Finally,
\begin{enumerate}
\item[(3)] For all $\mm \in \MaxSpec C$, $1 \leq \dim_\kk A/\mm A \leq d$,
\end{enumerate}
see e.g. \cite[Proposition III.1.1(2)]{BG-book}.

Denote the {\em{central character map}}
\begin{equation}
\label{kappa}
\kappa : \Irr(A) \to \MaxSpec C, \quad \kappa(V):= \Ann_C(V), \; \; \forall \, V \in \Irr(A).
\end{equation}
Its fibers are
\[
\kappa^{-1}(\mm) = \Irr(A/ \mm A), \; \; \forall \, \mm \in \MaxSpec A.
\]
\subsection{Cayley--Hamilton algebras}
\label{CH}
Denote by $\{\sigma_i \mid 1 \leq i \leq n \}$ the elementary symmetric functions in the variables
$\lambda_1, \lambda_2, \dots, \lambda_n$ and by
\[
\psi_i:=\lambda_1^i+\lambda_2^i+\cdots +\lambda_n^i
\]
the power sum symmetric functions.
As is well known, there exist unique polynomials
\[
g_i(x_1, x_2, \dots, x_i)\in \mathbb{Z}[(i!)^{-1}][x_1, x_2, \dots, x_i]
\]
such that
\[
\sigma_i=g_i(\psi_1, \psi_2, \dots, \psi_i), \quad \forall \; 1\leq i \leq n.
\]

Consider a unital algebra $A$ with trace over a field $\kk$ of characteristic $0$ or $>n$,
\[
\tr: A \to C,
\]
where $C$ is a subalgebra of $Z(A)$. Denote the $n$-th characteristic polynomial of an element $a \in A$:
\[
p_{n, a}(t):= t^n-c_1(a)t^{n-1}+\cdots +(-1)^n c_n(a) \in C[t],
\]
where $c_i(a):=g_i\big(\tr(a), \tr(a^2), \dots, \tr(a^i)\big)$.

\bde{CHalg} \cite{P,DP} Let $n$ be a positive integer $n$ and $\kk$ a field such that $\charr \kk \notin [1,n]$.
A Cayley--Hamilton algebra of degree $n$ is a unital $\kk$-algebra with trace $(A, C, \tr)$ such that
\begin{enumerate}
\item for all $a\in A$, $p_{n,a}(a)=0$ and
\item $\tr(1)=n$.
\end{enumerate}
\ede
\bex{CH-matr} (a) For every commutative unital algebra $C$, the standard trace function
\[
\tr : M_n(C) \to C
\]
makes the triple $(M_n(C), C, \tr)$ into a Cayley--Hamilton algebra of degree $n$ by the
Cayley--Hamilton identity for matrices with entries in a commutative algebra.

(b) Assume that $(A, C, \tr)$ is a Cayley--Hamilton algebra of degree $n$ and $A'$ is a unital subalgebra with a subalgebra
$C' \subseteq Z(A')$ such that $\tr(A') \subseteq C'$. It is straightforward to verify that
\[
(A', C', \tr|_{A'})
\]
is also a Cayley--Hamilton algebra of degree $n$.
%

(c) Assume that $A$ is a unital algebra and $C \subseteq Z(A)$ is a subalgebra such that $A$ is a free $C$-module of rank $n$.
Every $C$-basis $\{b_1, \ldots, b_n \}$ of $A$ gives rise to an algebra homomorphism
\[
\iota : A \to M_n(C), \quad \iota(a) := (c_{ij}) \in M_n(C),
\]
where $c_{ij} \in C$ are the unique elements with $a b_j = \sum_i c_{ij} b_i$.
Furthermore, $\ker \iota = 0$, since $a \in \ker \iota$ implies that $a b_i = 0$ for all $1 \leq i \leq n$,
and thus $a \big(\sum_i c_i b_i \big) =0$ for all $c_i \in C$, so $a . 1 = a =0$, because $A = C b_1 + \ldots + C b_n$.

The composition
\[
\tr_{\reg} := \tr \circ \iota : A \to C
\]
is called the regular trace of $A$. Clearly, it does not depend on the choice of $C$-basis of $A$.
By parts (a) and (b), the triple $(A, C, \tr_{\reg})$ is a Cayley--Hamilton algebra of degree $n$.
\eex

\bde{CHHopfalg} \cite{DPRR} A Cayley--Hamilton Hopf algebra is a Hopf algebra $H$ with a central Hopf subalgebra $C$ and a trace map
$\tr : H \to C$ such that $(H, C, \tr)$ is a Cayley--Hamilton algebra.
\ede
\subsection{Automorphisms and characters of Hopf algebras}
\label{sec:aut}
Let $H$ be a Hopf algebra over a field $\kk$. The characters of $H$ are precisely the elements of the group $G(H^\circ)$ of group like elements
of the restricted dual Hopf algebra $H^\circ$, see e.g. \cite[Example 9.1.4]{Mont}.
The identity element of this group is the counit $\ep : H \to \kk$. The inverse of a character
$\chi$ is the composition
\begin{equation}
\label{chi--1}
\chi^{-1}:=\chi S : H \to \kk,
\end{equation}
where $S$ is the antipode of $H$. The characters of $H$ are in bijection with the isomorphism classes of 1-dimensional modules of $H$. That is, we have
the canonical inclusion
\begin{equation}
\label{incl}
G(H^\circ) \subseteq \Irr (H).
\end{equation}

We will use Sweedler's notation for coproducts
\[
\Delta(h) = \textstyle \sum h_{(1)} \otimes h_{(2)}.
\]

Let $C$ be an arbitrary Hopf subalgebra of $H$. Its counit will be denoted by
\begin{equation}
\label{C-counit}
\ol{\ep} = \ep|_C : C \to \kk.
\end{equation}
Denote by
\[
\Aut_{\kk-\alg}(H,C)
\]
the group of $\kk$-algebra (not Hopf algebra) automorphisms $\theta$ of $H$ that preserve $C$, i.e., $\theta(C) \subseteq C$.

For $\chi \in G(H^\circ)$ and $h \in H$, define
\begin{equation}
\label{winding-aut}
W_l(\chi)(h) := \sum \chi(h_{(1)}) h_{(2)} \quad \mbox{and} \quad W_r(\chi)(h) := \sum \chi(h_{(2)}) h_{(1)}.
\end{equation}
Both maps are well known to be automorphisms of $H$, called the {\em{left}} and {\em{right winding}}
automorphism of $H$ associated to the character $\chi$. Since $C$ is a Hopf subalgebra of $H$, $\Delta (C) \subseteq C \otimes C$,
and thus, $W_l(\chi)(C) \subseteq C$ and $W_r(\chi)(C) \subseteq C$. Therefore,
\[
W_l(\chi), W_r(\chi) \in \Aut_{\kk-\alg}(H,C).
\]
It is well known that
\[
W_l : G(H^\circ) \to  \Aut_{\kk-\alg}(H,C) \quad \mbox{and} \quad
W_r : G(H^\circ) \to  \Aut_{\kk-\alg}(H,C)
\]
are an injective group antihomomorphism and homomorphism, respectively. The images
\[
W_l (G(H^\circ)) \quad \mbox{and} \quad W_r (G(H^\circ)),
\]
which are subgroups of $\Aut_{\kk-\alg}(H,C)$, are called the {\em{left}} and {\em{right winding automorphism
groups}} of $H$. For additional background we refer the reader to \cite[Section 2.5]{BZ}.
\sectionnew{Fiberwise actions on irreducible modules of Hopf algebras and Theorem A}
\label{sec:act}
In this section, we investigate an action of the isomorphism classes of the irreducible modules of the identity fiber algebra $H/ \mm_{\ol{\ep}} H$
on the isomorphism classes of the irreducible modules of the fiber algebras $H/ \mm H$ for Hopf algebras $H$ that are module finite over
a central Hopf subalgebra $C$. We describe the critical situation for its stabilizers and prove Theorem A from the introduction.
\subsection{A fiberwise action on irreducible modules}
\label{act}
In this section, we restrict to the Hopf case of the setting of Section \ref{Noether-alg}, i.e.,
work in the following setting:
\medskip

(FinHopf) {\em{$H$ is a finitely generated Hopf algebra over an algebraically closed field $\kk$
and $C$ is a central Hopf subalgebra such that $H$ is module finite over $C$}}.
\medskip

In \cite[Section 4.1]{Br} and \cite[Definition III.4.1]{BG-book}, a triple $(H, Z(H), C)$ as above, with the additional assumption that $H$ is a prime ring,
is called a PI Hopf triple. We will not need the last assumption in this paper. Following \cite[Section III.4.1]{BG-book}, note that if $C_1$ and $C_2$
are two subalgebras of $H$ satisfying (FinHopf), then the subalgebra of $H$ generated by $C_1$ and $C_2$ has the same property.
The maximal subalgebra of $H$ with this property is called the {\em{Hopf center}} of $H$ in \cite[Section III.4.1]{BG-book}.

By the Artin–Tate Lemma, $C$ is a finitely generated commutative Hopf algebra, and thus is the coordinate ring of an affine group scheme with identity element $\mm_{\ol{\ep}}$, recall \eqref{C-counit}.

By Skryabin's Theorem \cite{Sk}, the antipode of $H$ is bijective, and we can consider both the left and right duals of
a finite dimensional $H$-module $V$, which will be denoted by $V^*$ and ${}^* V$.
They are given by the dual vector space of $V$, equipped respectively with the actions
\[
\langle h. \xi, v \rangle : = \langle \xi, S(h) . v \rangle \quad \mbox{and} \quad  \langle h . \xi , v \rangle := \langle \xi, S^{-1}(h) . v \rangle
\]
for all $h \in V$ and $\xi$ in the dual vector space of $V$.

For a character $\vp \in G(C^\circ)$, denote the maximal ideal
\[
\mm_\vp := \Ker \vp \in \MaxSpec C.
\]
For $\mm \in \MaxSpec C$, $C/\mm \cong \kk$ by Hilbert's Nullstellensatz; that is, $\mm = \mm_\vp$ for some $\vp \in G(C^\circ)$.
For any two characters $\vp, \psi \in G(C^\circ)$ of $C$, we have
\[
(\vp \otimes \psi) \Delta(\mm_{\vp \psi}) = 0, \quad \mbox{and thus}, \quad
\Delta(\mm_{\vp \psi}) \subseteq \mm_\vp \otimes C + C \otimes \mm_\psi.
\]
Therefore, the coproduct $\Delta : H \to H \otimes H$ induces the algebra homomorphisms
\begin{equation}
\label{Deltaphipsi}
\Delta : H/\mm_{\vp \psi}H \to H/ \mm_\vp H \otimes H / \mm_\psi H.
\end{equation}
They turn $H/\mm_{\ol{\ep}} H$ into a finite dimensional Hopf algebra and $H/\mm H$ into finite dimensional left and right
comodule algebras over $H/\mm_{\ol{\ep}} H$:
\begin{equation}
\label{coact}
\Delta: H/\mm H \to H/\mm_{\ol{\ep}}H \otimes H/ \mm H, \quad \Delta: H/\mm H \to H/ \mm H \otimes H/\mm_{\ol{\ep}} H.
\end{equation}
In \cite[Section 4.1]{Br} and \cite[Section III.4.4]{BG-book},
the Hopf algebra $H/ \mm_{\ol{\ep}} H$ is called the {\em{restricted Hopf algebra}} associated to the pair $(H,C)$
in analogy with restricted universal enveloping algebras in positive characteristic which arise as special cases.
It follows from \eqref{coact} that
$\modd(H/\mm_{\ol{\ep}} H)$ is a tensor category and $\modd(H/\mm H)$ are left and right module categories over it, i.e.,
we have the bifunctors
\begin{align}
\label{l-tens}
&\modd(H/\mm_{\ol{\ep}} H) \times \modd(H/\mm H) \to \modd(H/\mm H) \quad \mbox{and} \\
\label{r-tens}
&\modd(H/\mm H) \times  \modd(H/\mm_{\ol{\ep}} H) \to \modd(H/\mm H),
\end{align}
induced by \eqref{coact}. Both bifunctors will be denoted by $\otimes$ as they are tensor products of $H$-modules
that factor through the quotients $H \to H/ \mm H$ for $\mm \in \MaxSpec C$. Furthermore, the homomorphisms \eqref{Deltaphipsi}
give rise to the bifunctors
\begin{equation}
\label{coDeltaphipsi}
\otimes : \modd(H/ \mm_\vp H) \times \modd(H / \mm_\psi H) \to \modd(H/\mm_{\vp \psi}H)
\end{equation}
for all characters $\vp, \psi \in G(C^\circ)$.

The next theorem addresses the action of the identity fiber $\Irr(H/\mm_{\ol{\ep}} H)$ of the central character map $\kappa : \Irr(H) \to \MaxSpec C$
from \eqref{kappa} (for $A:=H$) on the other fibers $\Irr(H/\mm H)$ of $\kappa$.
\bth{action-irr} Assume that $(H,C)$ is a pair of a Hopf algebra and a central subalgebra satisfying {\em{(}}FinHopf\,{\em{)}}.
Let $\mm \in \MaxSpec C$.
\begin{enumerate}
\item[(a)] For all $V, W \in \Irr(H/\mm H)$, there exist $M', M'' \in \Irr(H/\mm_{\ol{\ep}} H)$, such that $W$ is a quotient of $M' \otimes V$
and a submodule of $M'' \otimes V$.
\item[(b)] The group
\[
G_0:=G( (H/\mm_{\ol{\ep}}H)^\circ ) \subseteq \Irr (H/\mm_{\ol{\ep}}H)
\]
acts on $\Irr(H/\mm H)$ by $V \in \Irr(H/\mm H) \mt \chi \otimes V$ for $\chi \in G_0$.
\item[(c)] {\em{\cite[Proposition III.4.11]{BG-book}}} If, in addition the algebra $H/\mm_{\ol{\ep}} H$ is a basic, then the inclusion \eqref{incl} is an equality
\[
G_0= G( (H/\mm_{\ol{\ep}}H)^\circ ) = \Irr (H/\mm_{\ol{\ep}}H)
\]
and the action in part (b) is transitive.
\end{enumerate}
\eth
Recall that, for all finite dimensional $H$-modules $M,V$ and $W$, we have the adjunction isomorphisms
\begin{align}
&\Phi' : \Hom_H(M, W \otimes V^*) \stackrel{\cong}{\longrightarrow} \Hom_H(M \otimes V, W) \; \; \mbox{and}
\label{adjunction1}
\\
&\Phi'' : \Hom_H(W \otimes {}^* V, M) \stackrel{\cong}{\longrightarrow} \Hom_H(W, M \otimes V),
\label{adjunction2}
\end{align}
given by
\[
\Phi'(\eta')(m \otimes v) := (\id_W \otimes v) \eta'(m) \quad \mbox{and} \quad
(\Phi'')^{-1}(\eta'')(w \otimes \xi) := (\id_M \otimes \xi) \eta''(w)
\]
for all $\eta' \in \Hom_H(M, W \otimes V^*)$, $\eta'' \in \Hom_H(W, M \otimes V)$ and $m \in M$, $v \in V$, $w \in W$, $\xi \in {}^* V$.
\begin{proof}[Proof of \thref{action-irr}]
(a) Since $V, W \in \modd(H/\mm H)$, $V^* \in \modd(H / S^{-1}(\mm) H)$, and by \eqref{coDeltaphipsi},
$W \otimes V^* \in \modd(H/ \mm_{\ol{\ep}} H)$. Choose an irreducible module $M' \in \modd(H/ \mm_{\ol{\ep}}H)$ in the socle
of $W \otimes V^*$, see e.g. \cite[Section I.3]{ASS} for terminology.
It gives rise to a nontrivial homomorphism
\[
\eta' \in \Hom_H(M', W \otimes V^*) .
\]
The first adjunction formula \eqref{adjunction1} produces a nonzero homomorphism
\[
\Phi'(\eta') \in \Hom_H(M' \otimes V, W).
\]
Since $W$ is irreducible, this homomorphism is surjective. Therefore, $W$ is a quotient of $M' \otimes V$. This proves the
first statement of part (a) of the theorem.

Similarly, we have that $W \otimes  {}^* V \in \modd(H/ \mm_{\ol{\ep}})$. Let $M'' \in \modd(H/ \mm_{\ol{\ep}})$
be an irreducible $H$-module
in the top of $W \otimes  {}^* V$, see e.g. \cite[Section I.3]{ASS} for terminology.
We obtain a nontrivial homomorphism
\[
\eta'' \in \Hom_H(W \otimes {}^* V, M'').
\]
The second adjunction formula \eqref{adjunction2} gives a nonzero homomorphism
\[
\Phi''(\eta'') \in \Hom_H(W, M'' \otimes V).
\]
This homomorphism is injective because $W$ is irreducible. Therefore, $W$ is a submodule of $M \otimes V$. This proves the
second statement of part (a) of the theorem.

(b) For all $\chi \in G_0$ and $V \in \Irr(H/ \mm H)$, $\chi \otimes V$ is an irreducible $H$-module since for every $H$-submodule $W \subseteq \chi \otimes V$,
$\chi^{-1} \otimes W$ is a submodule of $V$, recall \eqref{chi--1}. By \eqref{l-tens}, $\chi \otimes V \in \modd(H/ \mm H)$, and thus, $\chi \otimes V \in \Irr(H/ \mm H)$.

(c) A proof of this part was given in \cite[Proposition III.4.11]{BG-book}. We include a proof based on part (a) for completeness.
If the algebra $H/\mm_{\ol{\ep}} H$ is a basic, then all of its irreducible modules are 1-dimensional. So,
\[
G( (H/\mm_{\ol{\ep}}H)^\circ ) = \Irr (H/\mm_{\ol{\ep}}H).
\]
Fix $V \in \Irr(H/ \mm H)$. By part (a), every $W \in \Irr(H / \mm H)$ is a submodule of $\chi \otimes V$ for some $\chi \in G_0$. Since $\chi \otimes V$ is
itself irreducible, this forces $W = \chi \otimes V$. Therefore, the $G_0$-action on $\Irr(H/ \mm H)$ is transitive.
\end{proof}
\bre{right-act} One can prove analogous statements to those in \thref{action-irr}(a)
for the tensor products $V \otimes M$ for $V \in \Irr(H/ \mm H)$ and
$M \in \Irr(H / \mm_{\ol{\ep}} H)$. This gives rise to a right action of $G_0$ on $\Irr(H/ \mm H)$
with the properties in \thref{action-irr}(b)--(c). We will not need  those facts in this paper and leave the details to the reader.
\ere
\bre{Frobenius} (1) Under the assumption (FinHopf), $H$ is a Frobenius extension of $C$, \cite[Corollary III.4.7]{BG-book} and, as a consequence, $H/ \mm H$
are (finite dimensional) Frobenius algebras for all $\mm \in \MaxSpec C$, \cite[Corollary III.4.7.2]{BG-book},
(see e.g. \cite[Section III.4.7]{BG-book} for definitions).

(2) In \cite{BG} general results relating the Azumaya locus of $H$ and the ramification locus of the canonical map
$\MaxSpec Z(H) \to \MaxSpec C$ were obtained under the assumption (FinHopf).
\ere
\subsection{Stabilizers}
\label{stab}
For the rest of this section we continue to work under the assumption that $(H,C)$ is a pair of a Hopf algebra and a central subalgebra satisfying (FinHopf).

For $V \in \Irr(H)$, denote the stabilizer of $V$ under the action in \thref{action-irr}(b):
\[
\Stab_{G_0}(V) = \{ \chi \in G_0 \mid \chi \otimes V \cong V \}.
\]

\ble{stab} \hfill
\begin{enumerate}
\item[(a)] For $V, W \in \Irr(H)$ and $\chi \in G_0$, the following are equivalent:
\begin{enumerate}
\item[(i)] $\chi \otimes V \cong W$;
\item[(ii)] $\chi$ is a submodule of $W \otimes V^*$;
\item[(iii)] $\chi$ is a quotient of $W \otimes {}^* V$.
\end{enumerate}
\item[(b)] For $V, W \in \Irr(H)$ and $\chi \in G_0$, the following are equivalent:
\begin{enumerate}
\item[(i)] $\chi \in \Stab_{G_0}(V)$;
\item[(ii)] $\chi$ is a submodule of $V \otimes V^*$;
\item[(iii)] $\chi$ is a quotient of $V \otimes {}^* V$.
\end{enumerate}
\end{enumerate}
\ele
The proof of part (a) is similar to the proofs of \thref{action-irr}(a)--(b), using the adjunction formulas \eqref{adjunction1}--\eqref{adjunction2},
and is left to the reader. The second part of the lemma follows from the first by letting $W=V$.

\bpr{stab} For all pairs $(H,C)$ of a Hopf algebra and a central subalgebra satisfying {\em{(}}FinHopf\,{\em{)}}, and a module $V \in \Irr(H)$,
the following hold:
\begin{enumerate}
\item[(a)] $\bigoplus_{\chi \in \Stab_{G_0}(V)} \chi$ is a submodule of $V \otimes V^*$ and a quotient of $V \otimes {}^* V$.
\item[(b)] $|\Stab_{G_0}(V)| \leq \dim(V)^2$.
\end{enumerate}
\epr
\begin{proof} (a) By \leref{stab}(b), $\chi$ is a submodule of $V \otimes V^*$ for all $\chi \in \Stab_{G_0}(V)$.
Therefore, $\bigoplus_{\chi \in \Stab_{G_0}(V)} \chi$ is a submodule of $V \otimes V^*$. The second statement
is proved in a similar way by using \leref{stab}(b) or by dualizing the first statement and using that $\chi^* = \chi S = \chi^{-1}$
for $\chi \in G_0$, so
\[
\Big( \textstyle \bigoplus_{\chi \in \Stab_{G_0}(V)} \chi \Big)^* \cong \bigoplus_{\chi \in \Stab_{G_0}(V)} \chi.
\]

Part (b) follows from part (a) since
\[
\dim \Big( \textstyle \bigoplus_{\chi \in \Stab_{G_0}(V)} \chi \Big) = | \Stab_{G_0}(V)|.
\]
\end{proof}
\subsection{Actions of twisted group algebras}
\label{twisted-grp-alg}
Recall that, for a group $G$ and a 2-cocycle
\[
\gamma : G \times G \to \kk^*,
\]
one defines the twisted group algebra
\[
\kk_\ga G := \textstyle  \bigoplus_{g \in G} \kk g
\]
with the product
\[
g . h = \gamma(g,h) gh,
\]
expended by $\kk$-bilinearity. For 2-cocycles $\beta$ and $\gamma$ that are cohomologous, the corresponding $\kk$-algebras
$\kk_\beta G$ and $\kk_\gamma G$ are isomorphic. A 2-cocycle $\gamma$ is normalized, if
\[
\gamma(1, g) = \gamma(g, 1) = 1, \quad \forall \, g \in G.
\]
For such a cocycle, $1 \in G$ is the identity element of $\kk_\gamma G$. Each 2-cocycle is cohomologous to a
normalized one.

Given a $G$-graded algebra
\[
A = \textstyle \bigoplus_{g \in G} A_g
\]
and a 2-cocycle $\gamma : G \times G \to \kk^*$, one defines the 2-cocycle twist $A_\gamma$, which is a $\kk$-algebra isomorphic to $A$
as a vector space, equipped with the new product
\[
a \bullet b := \gamma(g,h) a b, \quad \forall \, g, h \in G, \; a \in A_G, b \in A_h.
\]
The group algebra $\kk G$ is naturally $G$-graded by putting $\kk g$ in degree $g$ for all $g \in G$,
and the twist of it by a 2-cocycle $\gamma$ is canonically isomorphic to the corresponding twisted group algebra
\[
\kk_\gamma G \cong (\kk G)_\gamma.
\]

Returning to our setting (FinHopf), for a subgroup $G \subseteq G(H^\circ)$, consider the smash product
\[
H \# \kk G, \quad \mbox{where} \quad \chi h = W_l(\chi)^{-1}(h) \chi,
\quad \forall \chi \in G(H^\circ), \; h \in H,
\]
recall the definition \eqref{winding-aut} of the left winding automorphism $W_l(\chi) \in \Aut(H,C)$ and the fact that
$W_l : G(H^\circ) \to  \Aut_{\kk-\alg}(H,C)$ is a group antihomomorphism. Here and below, by abuse of notation, 
we denote $h:= h \# 1$ and $\chi:= 1 \# \chi$ for $h \in H$, $\chi \in G(H^\circ)$.

The algebra $H \# \kk G$ is $G$-graded
by putting $H \chi$ in degree $\chi$ for $\chi \in G$. For a 2-cocycle $\gamma : G \times G \to \kk^*$, we denote
the $\kk$-algebra
\[
H \# \kk_\gamma G := (H \# \kk G)_\gamma.
\]
It contains $H$ and $\kk_\gamma G$ as subalgebras and is isomorphic to $H \otimes \kk_\gamma G$ as an
$(H, \kk_\gamma G)$-bimodule. Furthermore,
\begin{equation}
\label{comm-rel}
\chi h = W_l(\chi)^{-1}(h) \chi, \quad \forall \chi \in G(H^\circ), \; h \in H.
\end{equation}

Consider an irreducible $H$-module $V$ and denote by
\begin{equation}
\label{pi-act}
\pi : H \to \End_\kk (V)
\end{equation}
the $H$-action on $V$. The homomorphism $\pi$ is surjective by Kaplansky's Theorem
\cite[Section I.13.3]{BG-book}.
For each $\chi \in \Stab_{G_0}(V)$, fix a linear operator
\[
L_\chi \in \End_\kk(V),
\]
giving the $H$-module isomorphism $V \stackrel{\cong}{\longrightarrow} \chi^{-1} \otimes V$,
where in the right hand side we use the canonical identification of $\kk$-vector spaces $\kk \otimes_\kk V \cong V$.
Since the module $V$ is irreducible, $L_\chi$ is uniquely defined up to rescaling.

We view $G_0 = G((H/\mm_{\ol{\ep}} H)^\circ)$ as a subgroup of $G(H^\circ)$ by composing
the projection $H \to H / \mm_{\ol{\ep}} H$ with the characters $\chi : H/ \mm_{\ol{\ep}} H \to \kk$.

\bth{twisted-rep} Assume that $(H,C)$ is a pair of a Hopf algebra and a central subalgebra satisfying {\em{(}}FinHopf\,{\em{)}}.
\begin{enumerate}
\item[(a)]
For each irreducible $H$-module $V$, there exists a uniquely defined 2-cocycle
\[
\gamma_V : \Stab_{G_0}(V) \times \Stab_{G_0}(V) \to \kk^*
\]
such that the map
\[
\chi \in \Stab_{G_0}(V) \mt L_\chi \in \End_\kk(V)
\]
defines an injective algebra homomorphism
\[
L : \kk_{\gamma_V} \Stab_{G_0}(V) \to \End_\kk(V).
\]
\item[(b)] The homomorphisms $L$ and $\pi$ {\em{(}}from
\eqref{pi-act}{\em{)}} turn $V$ into an \\
$H \# \kk_{\gamma_V} \Stab_{G_0}(V)$-module.
\end{enumerate}
\eth
In other words, part (a) states that $L: \Stab_{G_0} \to {\mathrm{GL}}(V)$ is a section
of a projective representation of the group $\Stab_{G_0}$ with Schur multiplier
$\gamma : G \times G \to \kk^*$, see \cite{Kar} for background on projective
group representations.
\begin{proof} (a) For each two characters $\chi, \theta \in G_0$, $L_\chi L_\theta \in \End_\kk(V)$ gives an isomorphism
\[
V \stackrel{\cong}{\longrightarrow} \theta^{-1} \otimes \chi^{-1} \otimes V.
\]
Therefore,
\begin{equation}
\label{L-prod}
L_\chi L_\theta = \gamma(\chi, \theta) L_{\chi \theta}
\end{equation}
for a (unique) scalar $\gamma(\chi, \theta) \in \kk^*$.

Consider the canonical identification $V \otimes V^* \cong \End_\kk(V)$ as $\kk$-vector spaces and the embedding
\[
\textstyle \bigoplus_{\chi \in \Stab_{G_0}(V)} \chi \hookrightarrow V \otimes V^*
\]
from \prref{stab}. Under this identification, $L_\chi \in \End_\kk(V)$ corresponds to picking up a nonzero vector in the image of $\chi \hookrightarrow V \otimes V^*$.
We note that this can be used as an alternative definition of the operators $L_\chi$.
The identification and \prref{stab}(a) imply that $L_\chi \in \End_\kk(V)$ are linearly
independent for $\chi \in \Stab_{G_0}(V)$. The associativity property of the algebra $ \End_\kk(V)$ now implies that
$\gamma : G \times G \to \kk^*$ is a 2-cocycle and the map
\[
\chi \mt L_\chi \quad
\mbox{defines a $\kk$-algebra embedding} \quad \kk_\ga \Stab_{G_0}(V) \hookrightarrow \End_\kk (V).
\]

(b) This homomorphism and the homomorphism $\pi : H \to \End_\kk(V)$ from \eqref{pi-act} combine to give
an $H \# \kk_\gamma \Stab_{G_0}(V)$-module structure on $V$ because
\[
L_\chi \pi(h) = \pi \big( W_l(\chi)^{-1}(h) \big) L_\chi, \quad \forall h \, \in H, \; \chi \in \Stab_{G_0}(V),
\]
which follows from the fact that $L_\chi \in \End_\kk(V)$ gives an $H$-module isomorphism
$V \stackrel{\cong}{\longrightarrow} \chi^{-1} \otimes V$.
\end{proof}

\bre{coh-class} It follows from \eqref{L-prod} that rescaling $L_\chi$ for $\chi \in \Stab_{G_0}(V)$
changes the cocycle $\gamma_V$ by a coboundary, so the cohomology class
\[
[\gamma_V] \in H^2(\Stab_{G_0}(V), \kk^*)
\]
is an invariant of $V$.
\ere
\subsection{Maximally stable irreducible modules}
\label{max-stab}
By \prref{stab}(b),
\[
|\Stab_{G_0}(V)| \leq  \dim(V)^2.
\]
\bde{Max-stab} An irreducible $H$-module $V$ will be called {\em{maximally stable}} if
\[
|\Stab_{G_0}(V)| = \dim(V)^2.
\]
\ede
For example, all 1-dimensional $H$-modules (i.e., all characters) are maximally stable.
\ble{inv-max-stable} For each $\chi \in G_0$, $V \in \Irr(H/ \mm H)$ is maximally stable
if and only if $\chi \otimes V$ is maximally stable.
\ele
\begin{proof} This follows from the facts that
\[
\Stab_{G_0}( \chi \otimes V) = \chi \, \Stab_{G_0}(V) \, \chi^{-1}
\]
and
\[
\dim(\chi \otimes V) = \dim (V).
\]
\end{proof}

\bpr{max-stab} For all pairs $(H,C)$ of a Hopf algebra and a central subalgebra satisfying {\em{(}}FinHopf\,{\em{)}}, and $V \in \Irr(H)$,
the following are equivalent:
\begin{enumerate}
\item[(i)] $V$ is a maximally stable module;
\item[(ii)] $V \otimes V^*$ is a direct sum of nonisomorphic 1-dimensional $H$-modules;
\item[(iii)] $V \otimes V^* \cong \bigoplus_{\chi \in \Stab_{G_0}(V)} \chi$.
\end{enumerate}
\epr
\begin{proof}
(i) $\Rightarrow$ (iii)  \prref{stab}(a) implies that $\bigoplus_{\chi \in \Stab_{G_0}(V)} \chi$ is a submodule of $V \otimes V^*$.
If $V$ is maximally stable, then the two modules have the same dimension, and thus, are equal.

(iii) $\Rightarrow$ (ii) This is obvious.

(ii) $\Rightarrow$ (i) Assume that
\[
V \otimes V^* \cong \textstyle \bigoplus_{\chi \in T} \chi
\]
for some subset $T \subseteq G_0$. \leref{stab}(b) implies that $\Stab_{G_0}(V) =T$. Therefore,
\[
|\Stab_{G_0}(V)| = |T| = \dim(V \otimes V^*) = \dim(V)^2,
\]
and hence, $V$ is maximally stable.
\end{proof}

\bth{max-stab} Let $(H,C)$ be pair of a Hopf algebra and a central subalgebra satisfying {\em{(}}FinHopf\,{\em{)}},
and $V$ be a maximally stable irreducible module of $H$. Then the primitive quotient
\[
H/ \Ann_H(V)
\]
is isomorphic to a twisted group algebra
\[
\kk_{\gamma_V} \Stab_{G_0}(V)
\]
for the 2-cocycle $\gamma_V : \Stab_{G_0}(V) \times \Stab_{G_0}(V) \to \kk^*$ from \thref{twisted-rep}.
Both algebras are isomorphic to $\End_\kk(V)$.
\eth
The proof of the theorem will be given after some analysis of the theorem and an example. 
\bre{simple} This theorem places a strong constraint on what Hopf algebras admit maximally stable irreducible modules of dimension $> 1$.

If $G$ is a finite group and $\charr \kk \nmid |G|$, then the group algebra $\kk G$
is semisimple by Maschke's Theorem. Such a group algebra is simple, if and only if $G$ has order 1 because the blocks of $\kk G$ correspond
to the irreducible $G$-modules, which in turn correspond to the conjugacy classes of $G$ (and the identity element $1 \in G$ forms a
single conjugacy class).

On the other hand, there are nontrivial examples of twisted group algebras that are simple as described by the next example.
\ere
\bex{ZlZl} For a positive integer $\ell$ consider the abelian group
\[
\Lambda := \Zset/\ell \times \Zset / \ell
\]
and denote its standard generators by $a$ and $b$ in multiplicative notation. Let $\kk$ be an
algebraically closed field of characteristic that does not divide $\ell$. Consider the 2-cocycle
\begin{equation}
\label{gamma}
\gamma: \Lambda \times \Lambda \to \kk^*, \quad
\gamma (a^i b^j, a^k b^m) = \ee^{jk}, \quad \forall \, i, j, k, m \in \Zset/ \ell,
\end{equation}
where $\ee \in \kk$ is a primitive $\ell$-th root of unity.
The corresponding twisted group algebra is isomorphic to a truncated quantum torus:
\[
\kk_\gamma \Lambda \cong \frac{\kk \langle a, b \rangle}{( a^\ell -1, b^\ell-1, ba - \ee ab )} \cdot
\]
There is an isomorphism $\kk_\gamma \Lambda \cong M_\ell(\kk)$, given by
\[
b \mt \diag(1, \ee, \ldots, \ee^{\ell-1}), \quad a \mt (g_{st})_{s,t=1}^\ell,
\; \; g_{st} :=
\begin{cases}
1, &\mbox{if} \; \; t \equiv s + 1 (\modd \ell)
\\
0, &\mbox{otherwise}.
\end{cases}
\]
This follows from noting that the above formulas define a surjective homomorphism $\kk_\gamma \Lambda \to M_\ell(\kk)$,
which has to be bijective since both algebras have dimension $\ell^2$. 
\eex
\begin{proof}[Proof of \thref{max-stab}]
The homomorphism $L : \kk_\gamma \Stab_{G_0}(V) \to \End_\kk(V)$ is injective.
By the assumption that $V$ is maximally stable,
\[
\dim \kk_\gamma \Stab_{G_0}(V) = \dim(V)^2 = \dim \End_\kk(V).
\]
Hence, the homomorphism $L$ is an isomorphism:
\[
L : \kk_\gamma \Stab_{G_0}(V) \stackrel{\cong}{\longrightarrow} \End_\kk(V).
\]
Since $V$ is an irreducible $H$-module, Kaplansky’s Theorem \cite[Section I.13.3]{BG-book} implies that the homomorphism
\eqref{pi-act} induces an isomorphism
\[
H/\Ann_H(V) \cong \End_\kk(V).
\]
The above two isomorphisms give the statement of the theorem.
\end{proof}

Theorem A in the introduction is the combination of \thref{action-irr}, Propositions \ref{pstab} and \ref{pmax-stab}, and \thref{max-stab}.
\sectionnew{The lowest level set of the square dimension function of a Hopf algebra and Theorems B and C}
\label{sec:Sd}
In this section we prove two results that characterize the lowest level set of the square dimension function of a Hopf
algebra $H$ which is module finite over a central Hopf subalgebra $C$. They are in turn used to prove Theorems B and C
from the introduction.
\subsection{The square dimension function of a Noether algebra}
Return to the setting of Section \ref{Noether-alg} of a finitely generated algebra $A$ over an algebraically closed field $\kk$ and
a central subalgebra $C$ such that $A$ is a finitely generated module over $C$ generated by $d$ elements.

\bde{sq-dim} Define the {\em{square dimension function}} of the pair $(A,C)$,
\begin{equation}
\label{Sd}
\Sd : \MaxSpec C \to \Zset \quad \mbox{by} \quad \Sd(\mm) := \sum_{V \in \Irr(A/\mm A)} \dim(V)^2.
\end{equation}
\ede
\noindent
In terms of the central character map $\kappa : \Irr(A) \to \MaxSpec C$ from \eqref{kappa},
the right hand side of the defining equation for $\Sd(\mm)$ ranges over the preimages $\kappa^{-1}(\mm)$.

By property (3) in Section \ref{Noether-alg}, the function $\Sd$ takes values in $[1,d^2]$.

We will be concerned with describing the lowest level set of the square dimension function of
$(A,C)$,
\[
\Sd^{-1}(\ell) \subseteq \MaxSpec C, \quad \mbox{where} \quad
\ell := \min \Sd(\MaxSpec C)
\]
in the Hopf case.

\subsection{First theorem for the lowest level set of the square dimension function of a Hopf algebra}
\label{sec:first-thm}
Through the rest of the section we will work in the following setting:
\medskip

(FinHopfBasic) {\em{$H$ is a finitely generated Hopf algebra over an algebraically closed field $\kk$ and $C$ is a central Hopf subalgebra such that
$H$ is module finite over $C$ and the algebra $H/\mm_{\ol{\ep}}H$ is basic.}}

\bth{first-lowest-sd-set} Let $(H,C)$ be a pair of a Hopf algebra and a central Hopf subalgebra satisfying {\em{(}}FinHopfBasic\,{\em{)}}.
The following hold:
\begin{enumerate}
\item[(a)] For all $\mm \in \MaxSpec C$,
\[
\Sd(\mm) = \frac{|G_0| \dim(V)^2}{|\Stab_{G_0}(V)|}
\]
for any $V \in \Irr(H/ \mm H)$, recall the definition \eqref{G0} of the group $G_0$.
\item[(b)] The lowest level set of the square dimension function of $(H,C)$ is
\[
\ell = |G_0|.
\]
\item[(c)] The following are equivalent for $\mm \in \MaxSpec(C)$:
\begin{enumerate}
\item[(i)] $\mm$ belongs to the lowest level set $\Sd^{-1}(|G_0|)$ of the square dimension function of $(H,C)$;
\item[(ii)] There exists $V \in \Irr(H/ \mm H)$ that is maximally stable;
\item[(iii)] All modules $V \in \Irr(H/ \mm H)$ are maximally stable.
\end{enumerate}
\end{enumerate}
\eth
\begin{proof} (a) Fix $V \in \Irr(H/ \mm H)$.
By \thref{action-irr}(c), the action of $G_0$ on $\Irr(H/\mm H)$ is transitive. Therefore,
$\Irr(H/ \mm H)$ consists of the isomorphism classes of the modules $\chi \otimes V$, where $\chi$ runs
over a set of representatives in $G_0$ of the cosets $G_0/\Stab_{G_0}(V)$.
For each of them, $\dim (\chi \otimes V) = \dim (V)$. Thus,
\[
\Sd(\mm) = \frac{|G_0|}{|\Stab_{G_0}(V)|} \dim(V)^2,
\]
which proves (a).

(b) We have
\[
\Sd(\mm_{\ol{\ep}}) = | \Irr(H/ \mm H)| = |G_0|.
\]
Part (a) of the theorem and \prref{stab}(b) imply that for all $\mm \in \MaxSpec C$,
\[
\Sd(\mm) = \frac{|G_0| \dim(V)^2}{|\Stab_{G_0}(V)|} \geq |G_0|,
\]
which proves (b).

(c) \leref{inv-max-stable} and the transitivity of the $G_0$-action on $\Irr(H/ \mm H)$ from \thref{action-irr}(c)
imply that  (ii) $\Leftrightarrow$ (iii).
By parts (a) and (b), $\mm \in \MaxSpec C$ belongs to the lowest level set $\Sd^{-1}(|G_0|)$ if and only if
\[
|G_0| = \Sd(\mm) = \frac{|G_0|}{|\Stab_{G_0}(V)|} \dim(V)^2,
\]
i.e., if and only if $|\Stab_{G_0}(V)| = \dim(V)^2$ for one $V \in \Irr(H/ \mm H)$, and thus, for any
$V \in \Irr(H/ \mm H)$. This proves that (i) $\Leftrightarrow$ (ii).
\end{proof}

\subsection{Second theorem for the lowest level set of the square dimension function of a Hopf algebra}
\label{sec:second-thm}
\bth{second-lowest-Sd-set}
Let $(H,C)$ be a pair of a Hopf algebra and a central Hopf subalgebra satisfying {\em{(}}FinHopfBasic\,{\em{)}}.
\begin{enumerate}
\item[(a)] For $\mm \in \MaxSpec C$ the following are equivalent:
\begin{enumerate}
\item[(i)] The algebra $H/\mm H$ has a 1-dimensional module;
\item[(ii)] The algebra $H/\mm H$ is basic;
\item[(iii)] $H/\mm H$ and $H/\mm_{\ol{\ep}} H$ are isomorphic as $\kk$-algebras;
\item[(iv)] $\mm$ belongs to the orbit of $\mm_{\ol{\ep}}$ under the group $\Aut_{\kk-\alg}(H,C)$
of all $\kk$-algebra automorphisms that preserve $C$;
\item[(v)] $\mm$ belongs to the orbit of $\mm_{\ol{\ep}}$ under the left winding automorphism group $W_l(G(H^\circ))$;
\item[(vi)] $\mm$ belongs to the orbit of $\mm_{\ol{\ep}}$ under the right winding automorphism group $W_r(G(H^\circ))$.
\end{enumerate}
\item[(b)] If $\mm \in \MaxSpec C$ satisfies any of the six equivalent conditions in part {\em{(}}a{\em{)}}, then
it belongs to the lowest level set $Sd^{-1}(|G_0|)$ of the square dimension function of $(H,C)$.
\item[(c)] If every maximally stable irreducible $H$-module is 1-dimensional, then a maximal ideal $\mm \in \MaxSpec C$
belongs to the lowest level set $Sd^{-1}(|G_0|)$ of the square dimension function of $(H,C)$
if and only if it satisfies any of the six equivalents conditions in part {\em{(}}a{\em{)}}.
\end{enumerate}
\eth
\begin{proof} (a) We prove the implications
\[
\mbox{(i)} \Rightarrow \mbox{(v)} \Rightarrow \mbox{(iv)} \Rightarrow \mbox{(iii)} \Rightarrow \mbox{(ii)} \Rightarrow \mbox{(i)}.
\]
Analogously to the first two implications, one shows that $\mbox{(i)} \Rightarrow \mbox{(vi)} \Rightarrow \mbox{(iv)}$.
\medskip

[(i) $\Rightarrow$ (v)] Consider a 1-dimensional module of $H/\mm H$ and the character of $H/\mm H$ associated to it.
Denote by $\vp : H \to \kk$ its canonical lift to a character of $H$. Let
\[
\ol{\vp} := \vp|_C : C \to \kk.
\]
We have $\mm_{\ol{\vp}} = \Ker \ol{\vp} = \mm$ and
\[
\ol{\ep} \big( W_l (\vp) (c) \big) = \textstyle \sum \ol{\vp}(c_{(1)}) \ol{\ep}(c_{(2)}) = \ol{\vp}(c), \quad \forall c \in C.
\]
Therefore,
\[
W_l(\vp) (\mm) = W_l(\vp) (\mm_{\ol{\vp}}) = \mm_{\ol{\ep}},
\]
and hence, $\mm \in W_l(G(H^\circ)) (\mm_{\ol{\ep}})$.
(We note that the characters of $H$ used in this argument are different from the characters
used earlier, because the latter come from characters of $H/ \mm H$ while the former from
characters of $H/ \mm_{\ol{\ep}} H$.)
\medskip

[(v) $\Rightarrow$ (iv)] This implication is obvious since $W_l(G(H^\circ))$ is a subgroup of the automorphism 
group $\Aut_{\kk-\alg}(H,C)$.
\medskip

[(iv) $\Rightarrow$ (iii)] This implication follows from the fact that, if $\mm$ and $\mm'$ are in the same orbit of $\Aut_{\kk-\alg}(H,C)$, then
\[
H/ \mm H \cong H / \mm' H
\]
as $\kk$-algebras.
\medskip

[(iii) $\Rightarrow$ (ii)] Since $H/ \mm H \cong H / \mm_{\ol{\ep}} H$ and $H / \mm_{\ol{\ep}} H$ is a
basic algebra, the $\kk$-algebra $H/ \mm H$ is also basic.
\medskip

[(ii) $\Rightarrow$ (i)] This implication is obvious.

(b) If $\mm \in \MaxSpec C$ satisfies any of the six equivalent conditions in part (a) of the theorem, then $H/\mm H \cong H/\mm_{\ol{\ep}} H$, and thus,
\[
\Sd(\mm) = \Sd(\mm_{\ol{\ep}}) = |G_0|.
\]
By \thref{first-lowest-sd-set}(b), $\mm$ belongs to the lowest level set $\Sd^{-1}(|G_0|)$ of the square dimension function of $(H,C)$

(c) This part follows from \thref{first-lowest-sd-set}(c) by comparing condition (ii) in \thref{first-lowest-sd-set}(c)
and condition (i) in \thref{second-lowest-Sd-set}(a).
\end{proof}
\subsection{Proofs of Theorems B(a)-(c) and C}
\label{Theorems-B-C}
Assume that $(H,C,\tr)$ is a Cayley--Hamilton Hopf algebras such that $H$ is a finitely generated $\kk$-algebra and
the identity fiber $H/ \mm_{\ol{\ep}} H$ is a basic algebra. Theorem 4.5 in \cite{DP} implies that $H$ is module finite over $C$ and $C$ is
a finitely generated $\kk$-algebra. Therefore the assumptions in Theorems \ref{tfirst-lowest-sd-set} and \ref{tsecond-lowest-Sd-set} are satisfied.

Theorem 4.1(b) in \cite{BY} implies that
\begin{align}
\label{V-discr}
\VV(D_k(H/C, \tr)) &= \VV(MD_k(H/C, \tr))
\\
&= \big\{ \mm \in \MaxSpec C \mid \Sd(\mm) < k \big\},
\nn
\end{align}
recall \eqref{Sd}.
Theorems B(a)-(c) and C now follow from this relation and Theorems \ref{tfirst-lowest-sd-set} and \ref{tsecond-lowest-Sd-set}.
\subsection{Proof of Theorem B(d)}
\label{ThmBd}
For the convenience of the reader, we recall several notions.
\bde{AS-Gorenstein} \hfill
\begin{enumerate}
\item[(a)] An algebra $A$ is called {\em{left AS-Gorenstein}} if
\begin{enumerate}
\item[(i)] ${}_A A$ has a finite injective dimension, to be denoted by $d$ and
\item[(ii)] for every irreducible $A$-module $S$, $\Ext^i_A(S,A)=0$ for all $i \neq d$ and $\Ext^d_A(S,A)$
is an irreducible right $A$-module.
\end{enumerate}

Here and below${}_A A$ and $A_A$ denote the left and right $A$ modules, given by the left and right actions of $A$ on itself.

Right AS-Gorenstein algebras are defined analogously.
\item[(b)] An algebra $A$ is called {\em{AS-Gorenstein}} if it is both left  and right AS-Gorenstein.
\end{enumerate}
\ede
Every Noetherian, finitely generated, PI Hopf algebra $A$ over a field $\kk$ is AS-Gorenstein by
\cite[Theorems 0.1 and 0.2(1)]{WZ}. The assumptions in Theorems A-C imply that the Hopf algebras
considered in those theorems are AS-Gorenstein.
\bde{hom-int} Let $H$ be an AS-Gorenstein Hopf algebra of injective dimension $d$.
\begin{enumerate}
\item[(a)] The left and right homological integrals of $H$ are the bimodules
\[
\textstyle \int^l = \Ext^d_H(\ep, {}_H H) \quad \mbox{and} \quad \textstyle \int^r = \Ext^d_{H^{\opp}} (\ep, H_H),
\]
where in the first and second formula $\ep$ refers to the left and right $H$ modules given by the counit of $H$, respectively.
In the second formula $H^{\opp}$ refers to the opposite algebra.
\item[(b)] The integral order of $H$, denoted $io(H)$, is the minimal positive integer $n$ such that $(\int^r)^{\otimes n}$ is
isomorphic to the trivial $H$-bimodule $\ep$.
\end{enumerate}
\ede
\noindent
{\em{Proof of Theorem B(d)}}. Proposition 4.5(c) in \cite{BZ} implies that
\[
\textstyle \int^l \cong {}^1 \kk^{w_l(\chi)}
\]
for some character $\chi$ of $H$, where the right hand side denotes the trivial $H$-bimodule
given by the counit $\ep$, twisted on the right side by the winding automorphism $W_l(\chi)$.
By \cite[Theorem 0.3]{BZ}, the Nakayama automorphism $\nu$ of $H$ equals
\[
\nu = S^2 W_l(\chi).
\]
We now apply \cite[Proposition 4.4(b)]{BZ}, stating that $\nu$ acts trivially on $Z(H)$,
to deduce that $W_l(\chi)$ acts trivially on $C$. This is equivalent to saying that
\[
\chi|_C = \ep|_C.
\]
Therefore, $\chi \in G_0 = G((H/ \mm_{\ol{\ep}})^\circ)$.

If $H$ has a bijective antipode, then
\[
S(\textstyle \int^l) = \textstyle \int^r \quad \mbox{and} \quad S(\textstyle \int^r) = \textstyle \int^l
\]
by \cite[Lemma 2.1]{LWZ}.
Since we work in the PI situation, this is the case by \cite{Sk}, and thus, $io(H)$ equals the
minimal positive integer $n$ such that $(\int^l)^{\otimes n}$ is
isomorphic to the trivial $H$-bimodule $\ep$.

Hence, $io(H)$ equals the order of the element $\chi \in G_0$, which implies the first statement in Theorem B(d).

For the second statement in Theorem B(d), take $V \in \Irr(H/ \mm H)$ and apply Theorem A(b)-(c) to obtain
\[
\Sd(\mm) \geq \frac{|G_0|}{\Stab_{G_0}(V)} \dim(V)^2 \geq |G_0|.
\]
Equation \eqref{V-discr} now implies the second statement in Theorem B(d).
\qed
\sectionnew{Applications}
\label{sec:examples}
In this section we present applications of Theorems A--C from the introduction
to group algebras of central extensions of abelian groups and (big) quantum groups at roots of unity.
This illustrates the theorems and answers naturally arising questions from them.
\subsection{Finite dimensional group algebras of central extensions}
\label{subsec:fin-group}

Fix a positive integer $\ell$. Consider the abelian group
\[
\Lambda:= \Zset/\ell \times \Zset / \ell
\]
as in \exref{ZlZl}, written in multiplicative notation
\[
\Lambda := \langle a, b \mid a^\ell = b^\ell = 1, ab = ba \rangle,
\]
and the cyclic group
\[
\Delta := \Zset / \ell,
\]
also written in multiplicative notation
\[
\Delta:=  \langle c \mid c^\ell = 1 \rangle.
\]
It is easy to verify that the map
\begin{equation}
\label{beta}
\beta : \Lambda \times \Lambda \to \Delta, \quad \mbox{given by} \quad
\beta (a^i b^j, a^k b^m) = c^{jk}, \quad \forall \, i, j, k, m \in \Zset/ \ell
\end{equation}
is a 2-cocycle. Denote the corresponding central extension
\[
\Sigma : = \Delta  \rtimes_\beta \Lambda, 
\]
i.e.,
\[
\Sigma:= \{c^s a^i b^j \mid a^i b^j \in \Lambda, c^s \in \Delta \}
\]
with the product
\[
(c^s a^i b^j) . (c^t a^k b^m) = c^{s+t + jk} a^{i+k} b^{j+m}, \quad \forall i,j,k,m,s,t \in \Zset/ \ell.
\]

Let $\kk$ be an algebraically closed field of characteristic $0$ or $> \ell^2$. Consider the group algebras
\begin{align*}
H&:= \kk \Sigma \cong \frac{\kk \langle a, b, c \rangle}{(ac - ca, bc - cb, ba - c ab, a^\ell -1, b^\ell -1, c^\ell -1)} \quad \mbox{and}
\\
C&:= \kk \Delta \cong \frac{\kk \langle c \rangle}{(c^\ell -1)} \cdot
\end{align*}
Since $\Sigma$ is a central extension of $\Lambda$ by $\Delta$, $C \subseteq Z(H)$. It is easy to verify from the presentation of $H$
that
\[
C = Z(H).
\]

The Hopf algebra $H$ is a free module over the Hopf subalgebra $C$ of rank $\ell^2$ with basis
\[
\{ a^i b^j \mid i, j \in \Zset/ \ell \}.
\]
By \exref{CH-matr}(c), the regular trace function $\tr_{\reg} : H \to C$ makes the triple $(H,C, \tr_{\reg})$
a Cayley--Hamilton algebra of degree $\ell^2$.

Denote by $\ee \in \kk$ a primitive $\ell$-th root of 1. We have the isomorphism
\[
C \cong \kk[c]/(c^\ell -1),
\]
and thus,
\[
\MaxSpec C = \{ \mm_s \mid s \in \Zset/ \ell \}, \quad
\mbox{where} \quad \mm_s := (c - \ee^s).
\]
The kernel of the counit $\ol{\ep} : C \to \kk$ is $\mm_{\ol{\ep}} = \mm_0$.
For $s \in \Zset / \ell$ consider the 2-cocycles
\[
\beta_s: \Lambda \times \Lambda \to \kk^*, \quad
\beta_s (a^i b^j, a^k b^m) = \ee^{jks}, \quad \forall \, i, j, k, m \in \Zset/ \ell.
\]
For $s=1$ we recover the cocycle from \eqref{gamma}, $\beta_1 = \gamma$, and for $s =0$ the trivial one.

The fiber algebras of $H$ with respect to all maximal ideals $\mm_s \in \MaxSpec C$ are
truncated quantum tori, which can be viewed as twisted group algebras of $\Lambda$,
\[
H / \mm_s H =  \frac{\kk \langle a, b \rangle}{(a^\ell -1, b^\ell -1, ba - \ee^s ab)} \cong \kk_{\beta_s} \Lambda,
\quad \forall s \in \Zset/ \ell.
\]
In particular, the identity fiber algebra is isomorphic to the group algebra of the abelian group $\Lambda$,
\[
H/ \mm_{\ol{\ep}} H = H / \mm_0 H = \kk[a, b]/(a^\ell -1, b^\ell -1) \cong \kk \Lambda,
\]
and so is a basic algebra. The corresponding group $G_0$ from Theorem A(a) is the dual group $\widehat{\Lambda}$:
\[
G_0 = G( (H / \mm_0 H)^\circ) = \widehat{\Lambda} = \{ \chi_{i, j} \mid i,j \in \Zset/ \ell\},
\quad \chi_{i,j}(a) = \ee^i, \chi_{i,j}(b) = \ee^j.
\]
For $s \in \Zset /\ell$, denote
\[
d := \gcd (\ell, s).
\]
It is easy to verify that the map
\begin{equation}
\label{rep1}
b \mt \diag(1, \ee^s, \ldots, \ee^{(\ell/d-1)s}), \quad a \mt (g_{rt})_{r,t=1}^{\ell/d},
\; \; g_{rt} :=
\begin{cases}
1, &\mbox{if} \; \; t \equiv r + 1 (\modd \ell/d)
\\
0, &\mbox{otherwise}
\end{cases}
\end{equation}
is an irreducible $H / \mm_s H$-module of dimension $\ell/d$. It will be denoted by $V_s$. We have that 
\begin{equation}
\label{tens-isom}
\chi_{i,j} \otimes V_s \cong V_s, \quad \forall i, j \in \Zset/ \ell \; \; \mbox{such that} \; \; d | i, j.
\end{equation}
Indeed, the representation $\chi_{i,j} \otimes V_s$ can be identified with the action of  $H / \mm_s H$ on $\kk^{\ell/d}$ 
given by 
\begin{equation}
\label{rep2}
b \mt \diag(\ee^{j}, \ee^{j+s}, \ldots, \ee^{j+ (\ell/d-1)s}), \quad a \mt (\ee^i g_{rt})_{r,t=1}^{\ell/d}.
\end{equation}
Denote the standard basis of $\kk^{\ell/d}$ by $\{e_1, \ldots, e_{\ell/d} \}$. We will use the notation $e_t$ for $t \in \Zset$
where $t$ is taken modulo $\ell/d$. We have $d = \alpha s + \beta \ell$ for some $\alpha, \beta \in \Zset$. 
The map
\[
e_t \mt \ee^{-i t} e_{t + \alpha j/d} 
\]
for $1 \leq t \leq \ell/d$ 
provides an isomorphism from the representation \eqref{rep2} to the representation \eqref{rep1}.

The number of the elements $\chi_{i,j}$ of $G_0$ in \eqref{tens-isom} equals $\ell^2/d^2$. \prref{stab}(b) implies that $|\Stab_{G_0}(V_s)| \leq \dim(V_s)^2 = \ell^2/d^2$, and so, 
\[
\Stab_{G_0}(V_s) = \widehat{\Lambda}^d = \{ \chi_{i,j} \mid i, j \in d(\Zset/ \ell) \}.
\]
(Here and below, we use $d$-th powers because, by convention, all of our abelian groups are written in the multiplicative notation.)
The corresponding 2-cocycle from \thref{twisted-rep} is
\[
\gamma_{V_s} = \beta_s|_{\widehat{\Lambda}^d}
\]
in terms of the isomorphism $\widehat{\Lambda} \cong (\Zset / \ell) \times (\Zset / \ell) \cong \Lambda$. \thref{action-irr}(c) implies
\begin{equation}
\label{Irr-Hs}
\Irr(H / \mm H) = \{ \chi_{i,j} \otimes V_s \mid
i, j \; \mbox{running over a set of representatives of $\widehat{\Lambda}/ \widehat{\Lambda}^d$} \}.
\end{equation}
Hence, the square dimension function \eqref{Sd} is given by
\[
\Sd(\mm_s) = d^2 \frac{\ell^2}{d^2} = \ell^2, \quad \forall s \in \Zset/ \ell.
\]

Applying \eqref{V-discr} yields
\[
\VV(D_k(A/C, \tr)) = \VV(MD_k(A/C, \tr))
=
\begin{cases}
\varnothing, & k \leq \ell^2
\\
\MaxSpec C, & k > \ell^2.
\end{cases}
\]
Another consequence of \eqref{Irr-Hs} is that the group of characters of $H$ is
\[
G(H^\circ) = G((H/ \mm_0)^\circ) = G_0.
\]
\thref{second-lowest-Sd-set}(a) implies that the orbits of $\mm_0$ under the left winding automorphisms of $H$, the right winding automorphisms of $H$ and
$\Aut(H,C)$ consist of the maximal ideal $\mm_0$ alone.  Thus, we obtain the following:
\bth{f-d-grp-als} Let $\ell$ be a positive integer and $\kk$ be an algebraically closed field of characteristic 0 or $> \ell^2$.
\begin{enumerate}
\item[(a)] The Hopf algebra
\[
H:= \kk \; \Big( (\Zset/ \ell) \rtimes_\beta  ( \Zset/\ell \times \Zset / \ell ) \Big)
\]
is a free module over its center
\[
Z(H) = \kk (\Zset/ \ell)
\]
of rank $\ell^2$. This gives rise to the Cayley--Hamilton Hopf algebra
structure $(H, Z(H), \tr_{\reg})$ of degree $\ell^2$.
\item[(b)] $H$ has maximally stable representations of dimension $g$ for all positive divisors $g$ of $\ell$.
All irreducible modules of $H$ are maximally stable.
\item[(c)] The zero sets of the discriminant/modified discriminant ideals of $(H, Z(H), \tr)$ are
\[
\VV(D_k(H/Z(H), \tr)) = \VV(MD_k(H/Z(H), \tr))
=
\begin{cases}
\varnothing, & k \leq \ell^2
\\
\Omega_\ell, & k > \ell^2,
\end{cases}
\]
where
\begin{equation}
\Omega_\ell := \{ \omega \in \kk \mid \omega^\ell = 1\}.
\label{Omega}
\end{equation}
\item[(d)] The left and right winding automorphism groups of $H$, as well as the group $\Aut(H,Z(H))$, 
fix the identity element of the group $\MaxSpec Z(H) \cong \Omega_\ell$.
\end{enumerate}
\eth

Part (b) of the proposition shows that, in general, the maximally stable representations $V$ of a Hopf algebra $H$ can have arbitrary dimension $\dim(V)$
such that
\[
\dim(V)^2 \mid |G_0|.
\]

Part (d) of the proposition shows that the inclusion in Theorem C(b) can be proper. Even worse, it can happen that
\[
\Aut(H,C) . \mm_{\ol{\ep}} = \{ \mm_{\ol{\ep}} \},
\]
while the zero set of the lowest discriminant ideal of $(H,C, \tr)$ is the whole $\MaxSpec C$.
\subsection{Infinite dimensional group algebras of central extensions} 
\label{subsec:inf-group}
The phenomena from the previous
subsection are not restricted to the case of finite dimensional Cayley--Hamilton Hopf algebras. Here we describe
how they appear in the infinite dimensional case.

Fix again a positive integer $\ell$ and an algebraically closed field $\kk$ of characteristic $\charr \kk \not \in [1, \ell^2]$,
and retain the notation from the previous subsection.
Consider the free abelian group
\[
\widetilde{\Lambda} := \langle x,  y \mid xy = yx \rangle \cong \Zset \times \Zset
\]
and the projection
\[
\pi : \wt{\Lambda} \to \Lambda,
\quad \pi(x) = a, \pi(y) = b.
\]
Denote the 2-cocycle
\[
\widetilde{\beta}: = \beta \circ \pi: \widetilde{\Lambda} \times \widetilde{\Lambda} \to \Delta
\]
and consider the central extension
\[
\widetilde{\Sigma} : = \Delta \rtimes_\beta \widetilde{\Lambda}.
\]
Its group algebra is
\[
\widetilde{H}:= \kk \widetilde{\Sigma} \cong \frac{\kk \langle x^{\pm 1}, y^{\pm 1}, c \rangle}{(xc - cx, yc - cy, y x - c xy,  c^\ell -1)} \cdot
\]
Denote by $\widetilde{\Delta}$ the subgroup of $\widetilde{\Sigma}$ generated by $x^\ell, y^\ell$ and $c$.
Clearly, $\widetilde{\Delta} \cong \Zset \times \Zset \times (\Zset / \ell)$. One verifies that
\begin{equation}
\label{ZHtilde}
Z(\widetilde{H}) = \kk \widetilde{\Delta} \cong \frac{\kk [ X^{\pm 1}, Y^{\pm 1}, c ]}{(c^\ell -1)}
\end{equation}
where the isomorphism is given by $x^\ell \mt X$, $y^\ell \mt Y$ and $c \mt c$. This implies in particular
that $Z(\widetilde{H})$ is a Hopf subalgebra of $\widetilde{H}$. The group $\Sigma$ is isomorphic to the factor group
of $\widetilde{\Sigma}$ by the free abelian subgroup generated by $x^\ell, y^\ell$ by sending
$x \mt a$, $y \mt b$, $c \mt c$, so
\[
H \cong \widetilde{H}/ (x^\ell-1, y^\ell -1).
\]

The Hopf algebra $\widetilde{H}$ is a free module over $Z(\widetilde{H})$ of rank $\ell^2$ with basis
\[
\{ x^i y^j \mid 0 \leq i, j \leq \ell -1\},
\]
which gives rise to the regular trace function $\tr_{\reg} : H \to C$, making $(H,C, \tr)$
a Cayley--Hamilton algebra of degree $\ell^2$ by \exref{CH-matr}(c).

It follows from \eqref{ZHtilde} that
\[
\MaxSpec Z(\widetilde{H}) = \{ \mm_{u, v, s} \mid u, v \in \kk^*, s \in \Zset/ \ell \}, \quad
\mbox{where} \quad \mm_{u, v, s} := (x^\ell - u, y^\ell -v, c - \ee^s).
\]
The kernel of the counit $\ol{\ep} : Z(\widetilde{H}) \to \kk$ is
\[
\mm_{\ol{\ep}} = \mm_{1,1,0}.
\]
The fiber algebras of $\widetilde{H}$ are
\begin{align*}
\widetilde{H} / \mm_{u, v, s} \widetilde{H} &= \frac{\kk \langle x, y \rangle}{(x^\ell -u, y^\ell - v, y x  - \ee^s x y)}
\\
& \cong \frac{\kk \langle a, b \rangle}{(a^\ell -1, b^\ell -1, ba - \ee^s ab)} \cong \kk_{\beta_s} \Lambda,
\quad \forall u, v \in \kk^*, s \in \Zset/ \ell,
\end{align*}
where the first isomorphism is given by $x \mt \sqrt[\ell]{u} a$, $y \mt \sqrt[\ell]{v} b$. From here on,
$\sqrt[\ell]{u}$ and $\sqrt[\ell]{v}$ will denote two fixed choices of $\ell$-th roots of $u$ and $v$ in $\kk$.
The identity fiber algebra is basic:
\[
\widetilde{H}/ \mm_{\ol{\ep}} \widetilde{H} = \widetilde{H} / \mm_{1,1,0} \widetilde{H} = \kk[x, y]/(x^\ell -1, y^\ell -1) \cong \kk \Lambda.
\]
The group $G_0$ from Theorem A(a) is the dual group $\widehat{\Lambda}$:
\[
G_0 = G( ( \widetilde{H} / \mm_{1,1,0} \widetilde{H})^\circ) = \widehat{\Lambda} = \{ \chi_{i, j} \mid i,j \in \Zset/ \ell\},
\quad \chi_{i,j}(x) = \ee^i, \chi_{i,j}(y) = \ee^j.
\]
In the notation $d := \gcd (\ell, s)$, the map
\[
y \mt \sqrt[\ell]{v} \diag(1, \ee^s, \ldots, \ee^{(\ell/d-1)s}), \; \;  x \mt (g_{rt})_{r,t=1}^{\ell/d},
\; g_{rt} :=
\begin{cases}
\sqrt[\ell]{u}, &\mbox{if} \; \; t \equiv r + 1 (\modd \ell/d)
\\
0, &\mbox{otherwise}
\end{cases}
\]
defines an irreducible $\widetilde{H} / \mm_{u,v,s} \widetilde{H}$-module of dimension $\ell/d$, to be denoted by $V_{u,v,s}$, and
\[
\Stab_{G_0}(V_{u,v,s}) = \widehat{\Lambda}^d = \{ \chi_{i,j} \mid i, j \in d(\Zset/ \ell) \}.
\]
In particular, $V_{u,v,s}$ is maximally stable.
The corresponding 2-cocycle from \thref{twisted-rep} is again computed to be $\gamma_V = \beta_s|_{\widehat{\Lambda}^d}$,
assuming the identification $\widehat{\Lambda} \cong \Zset \times \Zset \cong \Lambda$.
\thref{action-irr}(c) implies
\[
\Irr(\widetilde{H} / \mm_{u,v,s} \widetilde{H}) = \{ \chi_{i,j} \otimes V_s \mid
i, j \; \mbox{running over a set of representatives of $\widehat{\Lambda}/ \widehat{\Lambda}^d$} \},
\]
and thus,
\[
\Sd(\mm_{u,v,s}) = \ell^2, \quad \forall u, v \in \kk^*, s \in \Zset/ \ell.
\]

Equation \eqref{V-discr} and \thref{second-lowest-Sd-set}(a) imply the following:

\bth{f-d-grp-als-inf-dim} Let $\ell$ be a positive integer, $\kk$ be an algebraically closed field of characteristic $\charr \kk \notin[1, \ell^2]$ and
\[
\widetilde{H}:= \kk \; \Big( (\Zset/ \ell) \rtimes_{\widetilde{\beta}} ( \Zset \times \Zset ) \Big).
\]
\begin{enumerate}
\item[(a)] The Hopf algebra $\widetilde{H}$ is free over its center $Z(\widetilde{H})$, which is isomorphic to the coordinate
ring of the affine algebraic group $\mathbb{G}_m^2 \times \Omega_\ell$ as a Hopf algebra. This makes the
triple $(H, Z(H), \tr_{\reg})$ a Cayley--Hamilton Hopf algebra of degree $\ell^2$.
\item[(b)] All irreducible modules of $\widetilde{H}$ are maximally stable and those of the
fiber algebras $H/ \mm_{u, v, s} H$ have dimension $\ell^2/\gcd(s, \ell)^2$
for all $u, v \in \kk^*$, $s \in \Zset/ \ell$.
\item[(c)] We have
\[
\VV(D_k(\widetilde{H}/Z(\widetilde{H}), \widetilde{\tr})) = \VV(MD_k(\widetilde{H}/Z(\widetilde{H}), \widetilde{\tr}))
=
\begin{cases}
\varnothing, & k \leq \ell^2
\\
\mathbb{G}_m^2 \times \Omega_\ell, & k > \ell^2.
\end{cases}
\]
\item[(d)] The winding automorphism groups $W_l(\widetilde{H})$ and the automorphism group $W_r(\widetilde{H})$ and the group 
$\Aut(\widetilde{H},Z(\widetilde{H}))$ fix the identity element of $\MaxSpec Z(\widetilde{H}) \cong \mathbb{G}_m^2 \times \Omega_\ell$.
\end{enumerate}
\eth
One can easily show that the integral orders of the Hopf algebras in this and the previous subsection equal 1, i.e., they are unimodular. 
\subsection{Big quantum Borel subgroups at roots of unity}
Let $\g$ be a finite dimensional, complex, simple Lie algebra of rank $r$ and $G$ be the corresponding connected, simply connected, complex algebraic group. Denote by
$B^\pm$ a pair of opposite Borel subgroups and by $T:= B^+ \cap B^-$ the corresponding maximal torus of $G$. The coordinate rings of $G$ and $B^-$ will
be denoted by $\OO(G)$ and $\OO(B^-)$. The Cartan matrix of $\g$ will be denoted by $(a_{ij})$.
Let $d_1, \ldots, d_r$ be the relatively prime, positive integers that symmetrize the Cartan matrix $(a_{ij})$.
Denote by $\Phi_+$ and $P_+$ the sets of positive roots and dominant integral weights of $\g$, respectively.
Let $\{\al\spcheck_1, \ldots, \al\spcheck_r \}$ be the simple coroots of $\g$.

Let $W$ be the Weyl group of $\g$ and $l : W \to \{0, 1, \ldots\}$ be its length function. The reflection length function $s : W \to \{0, 1, \ldots \}$
is defined by setting $s(w)$ to be equal to the minimal length of a presentation of $w \in W$ as a product of reflections. It also equals
\[
s(w) = \dim \coker (\id - w)
\]
for the action of $w$ on $\Lie T$. The double Bruhat cells of $G$ are defined by
\[
G^{w_1, w_2} := B^+ w_1 B^+ \cap B^- w_2 B^- \quad \mbox{for} \quad w_1, w_2 \in W.
\]

Let $\ee \in \Cset$ be a root of unity of odd order $\ord(\ee)$, which is coprime to 3 if $\g$ is of type $G_2$. In other words,
$\ord(\ee)$ is odd and coprime to $d_i$ for $1 \leq i \leq r$. Denote by $\UU_\ee(\g)$ the De Concini--Kac \cite{DK}
(big) quantized universal enveloping algebra of $\g$ at the root of unity $\ee$. It is a non-cocommutative Hopf algebra with
Chevalley generators
\[
\{K_i^{\pm 1}, E_i, F_i \mid 1 \leq i \leq r\}.
\]

Denote by $\UU_\ee(\g)^{\geq}$ the quantum Borel subalgebra of $\UU_\ee(\g)$, \cite{DKP-solv,DP}, 
generated by the Chevalley generators $E_i, K_i^{\pm 1}$,
$1 \leq i \leq r$, subject to the relations
\begin{align*}
K_i K_j = K_j K_i, \; \;
K_i E_j K_i^{-1} = \ee^{d_i a_{ij}} E_j, \; \;
\sum_{k=0}^{1-a_{ij}} (-1)^k
\begin{bmatrix}
1 - a_{ij} \\
k
\end{bmatrix}_{\ee^{d_i}}
E_i^{1 - a_{ij} -k} E_j E_i^k =0,
\end{align*}
where in the last identity $i \neq j$ and
\[
\quad [n]_{\ee^d} := \frac{\ee^{nd} - \ee^{-nd}}{\ee^d - \ee^{-d}}, \quad
[n]_{\ee^d}! := [1]_{\ee^d} \cdots [n]_{\ee^d} \quad \mbox{and} \quad
\begin{bmatrix}
n \\
k
\end{bmatrix}_{\ee^d}
:= \frac{[n]_{\ee^d}!}{ [k]_{\ee^d}! [n-k]_{\ee^d}!}
\]
for the positive integer values of $k \leq n$ for which the denominators involved do not vanish.
$\UU_\ee(\g)^{\geq}$ is a Hopf subalgebra of $\UU_\ee(\g)$ with coproduct, antipode and counit given by
\begin{gather*}
\Delta(K_i) = K_i \otimes K_i, \quad
\Delta(E_i) = E_i \otimes 1 + K_i \otimes E_i,
\\
S(K_i) = K_i^{-1}, \quad S(E_i) = - K_i^{-1} E_i,
\\
\ep(K_i) =1, \quad \ep(E_i) = 0.
\end{gather*}

De Concini and Kac \cite{DK} introduced the central Hopf subalgebra
\[
C_\ee(\g)^{\geq} := \Cset[ K_i^{\pm \ord(\ee)}, E_\alpha^{\ord(\ee)} \mid 1 \leq i \leq r, \alpha \in \Phi_+]
\subset Z(\UU_\ee(\g)^{\geq}),
\]
where $\{ E_\alpha \mid \alpha \in \Phi_+ \}$ are the positive root vectors of  $\UU_\ee(\g)$, \cite[Corollaries 3.1 and 3.3]{DK},
and interpreted it as the image of a quantum Frobenius map. There is a canonical Hopf algebra isomorphism
\begin{equation}
\label{two-algs}
C_\ee(\g)^{\geq} \cong \OO(B^-).
\end{equation}
\cite[Theorem 14.1]{DP} deals with the image of the quantum Frobenius map of $\UU_\ee(\g)$ and the above is its restriction
to $C_\ee(\g)^{\geq}$.
We will identify the maximal ideals of the two algebras in \eqref{two-algs} via this isomorphism. The maximal ideal of the counit
$\mm_{\ol{\ep}} \in \MaxSpec C_\ee(\g)^{\geq}$ corresponds to the identity element $1 \in B^-$.

The algebra $\UU_\ee(\g)^{\geq}$ is a free module over $C_\ee(\g)^{\geq}$ of rank
$\ord(\ee)^{\dim B^-}$ by using PBW bases. This gives rise to the regular trace function
\[
\tr_{\reg} : \UU_\ee(\g)^{\geq} \to \CC_\ee(\g)^{\geq},
\]
turning $(\UU_\ee(\g)^{\geq}, \CC_\ee(\g)^{\geq}, \tr_{\reg})$ into a Cayley--Hamilton Hopf algebra of degree equal to $\ord(\ee)^{\dim B^-}$,
cf. \exref{CH-matr}(c).

The (Jacobson) radical of a finite dimensional algebra $A$ will be denoted by $\rad(A)$, see e.g. \cite[Chapter I]{ASS}.
\ble{char} \hfill
\begin{enumerate}
\item[(a)] The characters of $\UU_\ee(\g)^{\geq}$ are given by
\begin{equation}
\label{char-Ugeq}
\chi_t(K_i) = t_i, \; \; \chi_t(E_i) =0, \quad \forall \; 1 \leq i \leq r
\end{equation}
for $t:=(t_1, \ldots, t_r) \in (\Cset^*)^r$.
\item[(b)] We have
\begin{multline*}
\Big( \UU_\ee(\g)^{\geq} / \mm_{\ol{\ep}}  \, \UU_\ee(\g)^{\geq} \Big)
/ \rad \Big( \UU_\ee(\g)^{\geq} / \mm_{\ol{\ep}} \,  \UU_\ee(\g)^{\geq} \Big)
\\
\cong \Cset [K_1, \ldots, K_r]/ (K_i^{\ord(\ee)} -1 \mid 1 \leq i \leq r).
\end{multline*}
In particular, the identity fiber algebra $\UU_\ee(\g)^{\geq} / \mm_{\ol{\ep}}  \UU_\ee(\g)^{\geq}$ is basic and
\begin{equation}
\label{G0-qB}
G_0 = G \big( ( \UU_\ee(\g)^{\geq} / \mm_{\ol{\ep}}  \UU_\ee(\g)^{\geq})^\circ \big)
\cong (\Zset/ \ord(\ee))^{r},
\end{equation}
recall \eqref{G0}.
\end{enumerate}
\ele
\begin{proof} (a) If $\chi : \UU_\ee(\g)^{\geq} \to \Cset^*$ is a character, then
\[
\chi(K_i) = t_i, \quad \forall \; 1 \leq i \leq r
\]
for some $t_i \in \Cset^*$ because the elements $K_i$ are units. The assumption that the order of the root of unity $\ee$ is odd and coprime to the
symmetrizing integers $d_i$ of the Cartan matrix $a_{ij}$ implies that $ \ee^{d_i a_{ij}} \neq 1$. The
identity $K_i E_i K_i^{-1} = \ee^{2 d_i} E_i$ implies $\chi(E_i) = 0$.

It is obvious from the presentation of $ \UU_\ee(\g)^{\geq}$ that \eqref{char-Ugeq} defines a character of  $\UU_\ee(\g)^{\geq}$
for all $t_1, \ldots, t_r \in \Cset^*$.

(b) By the definition of the counit of $ \UU_\ee(\g)^{\geq}$,
\[
\mm_{\ol{\ep}} = ( K_i^{\ord(\ee)} -1, E_\alpha^{\ord(\ee)} \mid 1 \leq i \leq r, \alpha \in \Phi_+),
\]
and thus,
\[
E_\alpha \in \rad \big( \UU_\ee(\g)^{\geq} / \mm_{\ol{\ep}} \,  \UU_\ee(\g)^{\geq} \big).
\]
Therefore, we have a surjective homomorphism
\begin{multline*}
\Cset [K_1, \ldots, K_r]/ (K_i^{\ord(\ee)} -1 \mid 1 \leq i \leq r) \to
\\
\Big( \UU_\ee(\g)^{\geq} / \mm_{\ol{\ep}}  \, \UU_\ee(\g)^{\geq} \Big)
/ \rad \Big( \UU_\ee(\g)^{\geq} / \mm_{\ol{\ep}} \,  \UU_\ee(\g)^{\geq} \Big).
\end{multline*}
This homomorphism is an isomorphism since $\Cset [K_1, \ldots, K_r]/ (K_i^{\ord(\ee)} -1 \mid 1 \leq i \leq r)$ is a semisimple subalgebra of
$\UU_\ee(\g)^{\geq} / \mm_{\ol{\ep}}  \, \UU_\ee(\g)^{\geq}$. This proves the first statement in part (b) of the lemma.
The last two statements of part (b) follow at once from it.
\end{proof}
We will need the following result.
\bth{rep-B} \cite[Theorem 2.3(a)]{BG0} Let
\[
\mm \in G^{w, 1} = B^+ w B^+ \cap B^- \subseteq B^- \cong \MaxSpec C_\ee(\g)^{\geq}
\]
for  $w \in W$. There are precisely $\ord(\ee)^{r - s(w)}$ nonisomorphic irreducible modules of
the fiber algebra $\UU_\ee(\g)^{\geq} / \mm \,  \UU_\ee(\g)^{\geq}$ and they all
have dimension $\ord(\ee)^{(l(w) + s(w))/2}$.
\eth

The next theorem shows that the big quantum Borel subgroups at roots of unity $\UU_\ee(\g)^{\geq}$ satisfy the assumptions of Theorem C(c) of the introduction and
derives consequences of this fact for the zero sets of their discriminant ideals.
\bth{quant-Borel} Let $\g$ be a finite dimensional complex simple Lie algebra of rank $r$ and $\ee \in \Cset$ be a root of unity of odd order,
coprime to $3$ if $\g$ is of type $G_2$. The following hold:
\begin{enumerate}
\item[(a)] $(\UU_\ee(\g)^{\geq}, \CC_\ee(\g)^{\geq}, \tr_{\reg})$ is a finitely generated Cayley--Hamilton Hopf algebra of degree $\ord(\ee)^{\dim B^-}$ with the properties that
\begin{enumerate}
\item[(i)] the identity fiber $\UU_\ee(\g)^{\geq} / \mm_{\ol{\ep}} \,  \UU_\ee(\g)^{\geq}$ is a basic algebra and
\item[(ii)] all maximally stable irreducible modules of $\UU_\ee(\g)^{\geq}$ have dimension 1.
\end{enumerate}
\item[(b)] The lowest discriminant ideal of $(\UU_\ee(\g)^{\geq}, \CC_\ee(\g)^{\geq}, \tr_{\reg})$ is of level
\[
\ell := \ord(\ee)^r +1
\]
and
\begin{align*}
&\VV(D_\ell(\UU_\ee(\g)^{\geq}/ \CC_\ee(\g)^{\geq}, \tr_{\reg})) = \VV(MD_\ell(  \UU_\ee(\g)^{\geq}/ \CC_\ee(\g)^{\geq}, \tr_{\reg}) )
\\
& =
W_l \big( G( (\UU_\ee(\g)^{\geq})^\circ)  \big) \, . \, \mm_{\ol{\ep}} =  W_r \big( G( (\UU_\ee(\g)^{\geq})^\circ) \big)  \, . \, \mm_{\ol{\ep}}
\\
&=
\Aut \big( \UU_\ee(\g)^{\geq},  \CC_\ee(\g)^{\geq} \big)  \, . \, \mm_{\ol{\ep}} = T \subset B^- \cong \MaxSpec \CC_\ee(\g)^{\geq}.
\end{align*}
\item[(c)] For $k > \ord(\ee)^r +1$,
\begin{multline*}
\VV(D_k(\UU_\ee(\g)^{\geq}/ \CC_\ee(\g)^{\geq}, \tr_{\reg})) = \VV(MD_k(  \UU_\ee(\g)^{\geq}/ \CC_\ee(\g)^{\geq}, \tr_{\reg}) )
\\
= \bigsqcup_{w \in W, \ord(\ee)^{r+ l(w)} < k} G^{w, 1}.
\end{multline*}
\end{enumerate}
\eth
\begin{proof} (a) The algebra $\UU_\ee(\g)^{\geq}$ is clearly finitely generated.
We proved that the triple $(\UU_\ee(\g)^{\geq}, \CC_\ee(\g)^{\geq}, \tr_{\reg})$ is a Cayley--Hamilton Hopf algebra
of degree $\ord(\ee)^{\dim B^-}$. Its identity fiber $\UU_\ee(\g)^{\geq} / \mm_{\ol{\ep}} \,  \UU_\ee(\g)^{\geq}$ is a basic algebra
by \leref{char}(b).

Assume that $V$ is an irreducible module of $\UU_\ee(\g)^{\geq}$ of dimension $>1$. \thref{rep-B} implies that
\[
V \in \Irr \big( \UU_\ee(\g)^{\geq} / \mm \, \UU_\ee(\g)^{\geq} \big)
\]
for some $\mm \in G^{w, 1}$ with $w \in W$, $w \neq 1$. Theorems A(a) and \ref{trep-B} and the formula \eqref{G0-qB}
for the group $G_0$ imply that
\[
\dim(V)^2 = \ord(\ee)^{l(w) + s(w)} \quad \mbox{and} \quad
|\Stab_{G_0} (V)| = \ord(\ee)^{r} / \ord(\ee)^{r - s(w)} = \ord(\ee)^{s(w)}.
\]
Hence, $|\Stab_{G_0} (V)| < \dim(V)^2$ because $l(w) \geq 1$, and thus, $V$ is not maximally stable.
This proves property (ii).

(b) The fact that the level of the lowest discriminant ideal of $\UU_\ee(\g)^{\geq}$ equals $\ord(\ee)^r +1$ follows from
Theorem B(b) and Equation \eqref{G0-qB}. The coproduct formulas for the generators of $\UU_\ee(\g)^{\geq}$
imply that for all $c = (c_1, \ldots, c_r) \in (\Cset^*)^r$,
\[
W_l(\chi_c) \, . \, \mm_{\ol{\ep}}  =
W_r(\chi_c) \, . \, \mm_{\ol{\ep}} =
\big( K_i^{\ord(\ee)} -c_i^{-1}, E_\alpha^{\ord(\ee)} \mid 1 \leq i \leq r, \alpha \in \Phi_+ \big),
\]
recall \eqref{char-Ugeq}. The description of the character group $G( (\UU_\ee(\g)^{\geq})^\circ)$
in \leref{char}(a) implies that
\[
W_l \big( G( (\UU_\ee(\g)^{\geq})^\circ)  \big)  \, . \, \mm_{\ol{\ep}} = W_r \big( G( (\UU_\ee(\g)^{\geq})^\circ)  \big) \, . \, \mm_{\ol{\ep}} = T.
\]
The other equalities in part (b) now follow from Theorem C(c).

Part (c) follows by combining Equation \eqref{V-discr} and \thref{rep-B}.
\end{proof}
\subsection{Quantized function algebras at roots of unity} We retain the notation from the previous subsection. In \cite{DL} De Concini and Lyubashenko
constructed a quantized coordinate ring $\OO_\ee(G)_{\Qset[\ee]}$ by first defining a rational form $\OO_q(G)_{\Qset[q^{\pm 1}]}$ of the quantized coordinate ring of $G$
over $\Qset[q^{\pm 1}]$ (a Hopf algebra over $\Qset[q^{\pm 1}]$) and then specializing $q$ to $\ee$:
\[
\OO_\ee(G)_{\Qset[\ee]} := \OO_q(G)_{\Qset[q^{\pm 1}]}/ (\Phi_{\ord(\ee)}(q)),
\]
where $\Phi_m(q) \in \Zset[q]$ denotes the $m$-th cyclotomic polynomial. As in \cite{BG0}, we will work over
$\Cset$ and $\Cset[q^{\pm 1}]$, respectively:
\[
\OO_\ee(G) := \OO_\ee(G)_{\Qset[\ee]} \otimes _{\Qset[\ee]} \Cset \cong \OO_q(G) /( q - \ee), \quad \mbox{where} \quad
\OO_q(G):= \OO_q(G)_{\Qset[q^{\pm 1}]} \otimes _{\Qset} \Cset.
\]
So, we have $\OO_\ee(G) \cong \OO_q(G)/(q-\ee)$. The role of the rational form $\OO_q(G)_{\Qset[q^{\pm 1}]}$ is only auxiliary to define 
$\OO_q(G)$, an algebra over $\Cset[q^{\pm 1}]$ and the specialization $\OO_\ee(G)$, an algebra over $\Cset$. 

In \cite[Proposition 6.4]{DL} a canonical central Hopf subalgebra $\CC_\ee(G)$ of $\OO_\ee(G)$ was constructed and it was proved
that
\[
\CC_\ee(G) \cong \OO(G)
\]
as Hopf algebras. We will use this isomorphism to identify
\[
\MaxSpec \CC_\ee(G) \cong G
\]
as algebraic groups. Under it, $\mm_{\ol{\ep}} \in \MaxSpec \CC_\ee(G)$ corresponds to the identity element $1 \in G$.

Furthermore, in \cite[Theorem 7.2]{DL} it was proved that $\OO_\ee(G)$ is a projective $\CC_\ee(G)$-module of rank $\ord(\ee)^{\dim G}$.
Based on this theorem, an improvement was obtained in \cite[Proposition 2.2]{BG0} and \cite[Theorem on p. 1]{BGS}
that $\OO_\ee(G)$ is a free module over $\CC_\ee(G)$ of rank $\ord(\ee)^{\dim G}$. This gives rise to the regular trace function
\[
\tr_{\reg} : \OO_\ee(G) \to \CC_\ee(G),
\]
turning $(\OO_\ee(G), \CC_\ee(G), \tr_{\reg})$ into a Cayley--Hamilton Hopf algebra of degree $\ord(\ee)^{\dim G}$,
cf. \exref{CH-matr}(c). We will need the following:

\bth{rep-O} \cite[Theorem 2.3(b)]{BG0} Let
\[
\mm \in G^{w_1, w_2} = B^+ w_1 B^+ \cap B^- w_2 B^- \subseteq G \cong \MaxSpec \OO_\ee(G)
\]
for  $w_1, w_2 \in W$. There are precisely $\ord(\ee)^{r - s(w_2^{-1}w_1)}$ nonisomorphic irreducible modules of
the fiber algebra $\OO_\ee(G) / \mm_{\ol{\ep}} \,  \OO_\ee(G)$ and they all
have dimension
\[
\ord(\ee)^{(l(w_1) + l(w_2) + s(w_2^{-1} w_1))/2}.
\]
\eth

Denote by $\UU_q(\g)$ the quantized universal enveloping algebra of $\g$ defined over $\Cset(q)$ and by $\UU_q(\g)_{\Cset[q^{\pm 1}]}$
the Lusztig divided power form of the quantized universal enveloping algebra defined over $\Cset[q^{\pm 1}]$. For a dominant
integral weight $\la \in P_+$, denote by $L(\la)$ the irreducible $\UU_q(\g)$-module of highest weight $\la$,
by $L(\la)_{\Cset[q^{\pm 1}]}$ its Lusztig integral form (which is a $\UU_q(\g)_{\Cset[q^{\pm 1}]}$-module) and its
dual
\[
L(\la)\spcheck_{\Cset[q^{\pm 1}]}:=\{ \xi \in L(\la)^* \mid \langle \xi, v \rangle \in \Cset[q^{\pm 1}], \forall v \in L(\la)_{\Cset[q^{\pm 1}]} \},
\]
where $L(\la)^*$ is the dual $\UU_q(\g)$-module.

The Hopf algebra $\OO_q(G)$ is the Hopf subalgebra of the restricted dual $(\UU_q(\g))^\circ$, given by the sums of all
matrix coefficients of the form
\[
c^\la_{\xi, v} \in (\UU_q(\g))^\circ, \quad \langle c^\la_{\xi, v}, x \rangle := \langle \xi,  x . v \rangle, \quad \forall x \in \UU_q(\g)
\]
for $\la \in P_+$, $v \in L(\la)_{\Cset[q^{\pm 1}]}$, $\xi \in L(\la)\spcheck_{\Cset[q^{\pm 1}]}$.
The product structure of $\OO_q(G)$ is given by
\begin{equation}
\label{prod-O}
c^{\la_1}_{\xi_1, v_1} c^{\la_2}_{\xi_2, v_2} := c^{\la_1 + \la_2}_{\xi_1 \otimes \xi_2, v_1 \otimes v_2}
\end{equation}
in terms of the canonical embedding of $L(\la_1 + \la_2)$ as a direct summand of $L(\la_1) \otimes L(\la_2)$.
This implies that the Hopf algebra $\OO_q(G)$ has the characters
\[
\chi_{\wt{t}} : \OO_q(G) \to \Cset[q^{\pm 1}], \quad
\mbox{given by} \quad \chi_{\wt{t}}( c^\la_{\xi, v} ) := \langle \xi, v \rangle \prod_{I=1}^r \wt{t}_i^{\langle \la, \al_i\spcheck \rangle}
\]
for
\[
\wt{t}:= (\wt{t}_1, \ldots, \wt{t}_r) \in \big( \Cset[q^{\pm 1}]^\times \big)^r.
\]
Here, as usual, $R^\times$ denotes the group of units of a ring $R$. 
This in turn, shows that the Hopf algebra $\OO_\ee(G)$ has the characters
\[
\chi_t : \OO_\ee(G) \to \Cset, \quad
\mbox{given by} \quad \chi_t( c^\la_{\xi, v} ) := \theta \big( \langle \xi, v \rangle \big) \prod_{I=1}^r t_i^{\langle \la, \al_i\spcheck \rangle}
\]
for
\[
t:= (t_1, \ldots, t_r) \in (\Cset^*)^r,
\]
where
\[
\theta : \Cset[q^{\pm 1}] \to \Cset
\]
is the $\Cset$-algebra homomorphism, given by $q \mt \ee$. Here we use the same notation for the elements of
$\OO_q(G)$ and their images in $\OO_\ee(G) = \OO_q(G) / (q - \ee)$. \thref{rep-O} implies that $\chi_t$ exhaust all characters
of $\OO_\ee(G)$. This can be also shown in an elementary way using \eqref{prod-O} as in \leref{char}(a).

The next theorem shows that all maximally stable modules of the root of unity quantum function algebras $\OO_\ee(G)$ 
are 1-dimensional, and so, Theorem C(c) is applicable to them.

\bth{O-eps} Let $\g$ be a finite dimensional complex simple Lie algebra of rank $r$ and $\ee \in \Cset$ be a root of unity of odd order,
coprime to $3$ if $\g$ is of type $G_2$. The following hold:
\begin{enumerate}
\item[(a)] $(\OO_\ee(G), \CC_\ee(G), \tr_{\reg})$ is a finitely generated Cayley--Hamilton Hopf algebra of degree $\ord(\ee)^{\dim G}$ with the properties that
\begin{enumerate}
\item[(i)] the identity fiber $\OO_\ee(G) / \mm_{\ol{\ep}} \, \OO_\ee(G)$ is a basic algebra and
\item[(ii)] all maximally stable irreducible modules of $\OO_\ee(G)$ have dimension 1.
\end{enumerate}
\item[(b)] The lowest discriminant ideal of $(\OO_\ee(G), \CC_\ee(G), \tr_{\reg})$ is of level
\[
\ell := \ord(\ee)^r +1
\]
and
\begin{align*}
&\VV(D_\ell(\OO_\ee(G)/ \CC_\ee(G), \tr_{\reg})) = \VV(MD_\ell( \OO_\ee(G)/ \CC_\ee(G), \tr_{\reg}) )
\\
& =
W_l \big( G( \OO_\ee(G)^\circ ) \big)  \, . \, \mm_{\ol{\ep}} =  W_r \big( G(\OO_\ee(G)^\circ ) \big)  \, . \, \mm_{\ol{\ep}}
\\
&=
\Aut \big( \OO_\ee(G), \CC_\ee(G) \big)  \, . \, \mm_{\ol{\ep}} = T \subset G \cong \MaxSpec \CC_\ee(G).
\end{align*}
\item[(c)] For $k > \ord(\ee)^r +1$,
\begin{multline*}
\VV(D_k( \OO_\ee(G)/ \CC_\ee(G), \tr_{\reg})) = \VV(MD_k( \OO_\ee(G) / \CC_\ee(G), \tr_{\reg})
\\
= \bigsqcup_{{w_1, w_w} \in W \times W, \ord(\ee)^{r+ l(w_1) + l(w_2)} < k} G^{w_1, w_2}.
\end{multline*}
\end{enumerate}
\eth
\begin{proof} (a) It is well known that the algebra $\OO_\ee(G)$ is finitely generated.
We proved that $(\OO_\ee(G), \CC_\ee(G), \tr_{\reg})$ is a Cayley--Hamilton Hopf algebra
of degree $\ord(\ee)^{\dim G}$. \thref{rep-O} implies that its identity fiber $\OO_\ee(G) / \mm_{\ol{\ep}} \,  \OO_\ee(G)$
is a basic algebra; this can be also derived directly from the definitions of $\OO_\ee(G)$ and $\CC_\ee(G)$.

The classification of the characters of $\OO_\ee(G)$ that we obtained now implies that
\[
G_0 = G \big( ( \OO_\ee(G) / \mm_{\ol{\ep}}  \OO_\ee(G))^\circ \big) =
\{ \chi_t \mid t \in \Om_{\ord(\ee)}^r \} \cong (\Zset/ \ord(\ee))^{\times r},
\]
recall \eqref{Omega}.

Let $V$ be an irreducible $\OO_\ee(G)$ -module of dimension $>1$. \thref{rep-O} implies that
\[
V \in \Irr \big( \OO_\ee(G) / \mm \OO_\ee(G) \big)
\]
for some $\mm \in G^{w_1, w_2}$ with $w_1, w_2 \in W$, such that either $w_1 \neq 1$ or $w_2 \neq 1$, i.e.,
$l(w_1) + l(w_2) > 0$.
Applying Theorems A(a) and \ref{trep-B} and  the equality $|G_0|= \ord(\ee)^r$, we obtain
\begin{align*}
&\dim(V)^2 = \ord(\ee)^{l(w_1) + l(w_2)+ s(w_2^{-1}w_1)} >  \ord(\ee)^{s(w_2^{-1}w_1)}
\\
&= \ord(\ee)^{r} / \ord(\ee)^{r - s(w_2^{-1} w_1)} = |\Stab_{G_0} (V)|.
\end{align*}
Therefore, $V$ is not maximally stable, which proves property (ii).

The first statement in part (b) follows from Theorem B(b) and the fact that the group $G_0$ has order $\ord(\ee)^r$,
Next we look at the action of the group of winding automorphisms of  $\OO_\ee(G)$. The coproduct of $\OO_q(G)$
is given by
\[
\Delta ( c^\la_{\xi, v} ) = \textstyle \sum_j  c^\la_{\xi, v_j} \otimes c^\la_{\xi_j, v}, \quad \forall v \in L(\la), \xi \in L(\la)\spcheck
\]
where $\{v_j \}$ and $\{\xi_j \}$ are dual bases of $L(\la)$ and $L(\la)\spcheck$, respectively. Therefore, for
$v \in L(\la)$ and $\xi \in L(\la)\spcheck$ of weights $\mu$ and $-\nu$, respectively,
\[
W_l( \chi_{\wt{t}}) ( c^\la_{\xi, v} ) =
\begin{cases}
\prod_{i=1}^r \wt{t}_i^{\langle \mu, \al_i\spcheck \rangle}, & \mu = \nu
\\
0, & \mbox{otherwise}
\end{cases}
\]
for all $\wt{t}= (\wt{t}_1, \ldots, \wt{t}_r) \in \big( \Cset[q^{\pm 1}]^\times \big)^r$. This implies that in $\OO_\ee(G)$, we have
\[
W_l( \chi_t) ( c^\la_{\xi, v} ) =
\begin{cases}
\prod_{i=1}^r t_i^{\langle \mu, \al_i\spcheck \rangle}, & \mu = \nu
\\
0, & \mbox{otherwise}
\end{cases}
\]
for all $t= (t_1, \ldots, t_r) \in (\Cset^*)^r$. Hence,
\[
W_l \big( G( \OO_\ee(G)^\circ ) \big)  \, . \, \mm_{\ol{\ep}} = T \subset G.
\]
The other equalities now follow from Theorem C(c).

Part (c) follows from Equation \eqref{V-discr} and \thref{rep-O} because
\[
\Sd(\mm) = \ord(\ee)^{r - s(w_2^{-1}w_1)} \Big( \ord(\ee)^{(l(w_1) + l(w_2) + s(w_2^{-1} w_1))/2} \Big)^2 =
\ord(\ee)^{r + l(w_1) + l(w_2) }
\]
for all $\mm \in G^{w_1, w_2}$ and $w_1, w_2 \in W$.
\end{proof}

\end{document}